\documentclass{amsart}
\usepackage{graphicx}
\usepackage{amssymb}
\usepackage{amscd}
\usepackage{latexsym}
\usepackage{amsfonts}
\usepackage{ragged2e}
\usepackage{amsmath}
\usepackage{float}
\usepackage[latin1]{inputenc}
\usepackage[english]{babel}

\usepackage{amsthm}

\newtheorem{proposition}{Proposition}
\newtheorem{lemma}[proposition]{Lemma}
\newtheorem{theorem}[proposition]{Theorem}
\newtheorem{corollary}[proposition]{Corollary}
\newtheorem{claim}[proposition]{Claim}
\newtheorem{thm}[proposition]{Theorem}

\def\Ahh{\frak A}

    \def\Err{{\bf R}}  \def\Zet{{\bf Z}}
        
 \def\CK{{\mathcal K}}  \def\CM{{\mathcal M}} 
\let\eps=\varepsilon
\let\alp=\alpha
\let\lam=\lambda

\def\supp{\mathop{\rm supp}\nolimits}

\def\isto{\buildrel \sim \over \longrightarrow}

\def\aff{{\mathop{\rm aff}\nolimits}}

\def\iff{\; \Longleftrightarrow \;}

\def\Proof{\noindent {\it Proof\/}}
\def\Remark{\medbreak \noindent {\it Remark\/}}

\def\satz #1: #2\par {{\bf #1:\enspace}{\it #2} \par
\ifdim\lastskip<\medskipamount \removelastskip\penalty55\medskip\fi}
\def\Satz #1: #2\par {\medbreak
     \noindent {\bf #1:\enspace}{\it #2} \par
\ifdim\lastskip<\medskipamount \removelastskip\penalty55\medskip\fi}
\def\teil #1 \par{\noindent {\it #1}\par
\ifdim\lastskip<\medskipamount \removelastskip\penalty55\medskip\fi}
\def\Teil #1) #2\par {\medbreak
\noindent {#1)\enspace}{\it #2} \par
\ifdim\lastskip<\medskipamount \removelastskip\penalty55\medskip\fi}

\newcommand{\ra}{\mathop{\fam0 \rightarrow}\nolimits}

\begin{document}

\title[Identities for Poincar\'e polynomials via Kostant cascades]{Identities for Poincar\'e polynomials \\ via Kostant cascades}

\date{\today}
\author{Jørgen Ellegaard Andersen}
\address{Center for Quantum Geometry of Moduli
  Spaces\\ 
  Department of Mathematics\\
  University of Aarhus\\
  DK-8000, Denmark}
	\email{andersen@qgm.au.dk}
	
  \author{Jens Carsten Jantzen}
	\address{Department of Mathematics\\
  University of Aarhus\\
  DK-8000, Denmark}
	\email{jantzen@math.au.dk}

\author{Du Pei}
\address{Walter Burke Institute for Theoretical Physics\\
 California Institute of Technology\\
Pasadena, CA 91125, USA\\ \newline \phantom{bb}
Center for Quantum Geometry of Moduli
  Spaces\\ 
  Department of Mathematics\\
  University of Aarhus\\
  DK-8000, Denmark}
\email{pei@caltech.edu}

\begin{abstract}
We propose and prove an identity relating the Poincar\'e polynomials of stabilizer subgroups of the affine Weyl group and of the corresponding stabilizer subgroups of the Weyl group. 
\end{abstract}

\thanks{Supported in part by the center of excellence grant ``Center for Quantum Geometry of Moduli Space" from the Danish National Research Foundation (DNRF95), by the Walter Burke Institute for Theoretical Physics, and by the U.S.~Department of Energy, Office of Science, Office of High Energy Physics, under Award Number {DE}-{SC0011632}.}

\maketitle

\section{Introduction}

Let $\Phi$ be a finite indecomposable root system and $W$ be its Weyl group. Recall that the affine Weyl group $W^{\text{aff}}$ is the semidirect product of $W$ and the root lattice ${\mathbb Z}\Phi$. It acts on $\mathbf{R}\Phi$ as a reflection group with (closed) alcoves as fundamental domains. Denote by  $\frak{A}$ the alcove contained in the dominant Weyl chamber such that $0\in \frak{A}$.

For $\lambda\in \frak{A}$, we consider the Poincar\'e polynomials of both its stabiliser $W_\lambda$ in $W$
$$P(W_\lambda) = \sum_{w\in W_\lambda}t^{l(w)}$$
and its stabiliser $W_\lambda^{\text{aff}}$ in $W^{\text{aff}}$
$$P\left(W_\lambda^{\text{aff}}\right) = \sum_{w\in W_\lambda^{\text{aff}}}t^{l(w)},$$
where $l$ is the length function.
Now $W_\lam$ is a parabolic subgroup of $W_\lambda^{\text{aff}}$, so $P_\lambda(t)$ divides $P_\lambda^{\text{aff}}(t)$, and we will establish a decomposition of the quotient $P\left(W_\lambda^{\text{aff}}\right)/P(W_\lambda)$ in terms of stabilisers. 

For this purpose, define $W_{\lambda,\xi}:=\mathrm{Stab}_{W_\lambda}(\xi)$ in $W_\lambda$, where $\xi\in \mathbf{Z}\Phi$ is in the root lattice. Then 
$$
P\left(W_\lambda/W_{\lambda,\xi}\right):=P(W_\lambda)/P(W_{\lambda,\xi})
$$
is again a polynomial. We then define an integral function on $\mathbf{Z}\Phi$ given by
$$
\ell_\xi:=\mathrm{deg}\left(P(W_\xi^{\text{aff}})/P(W_{\xi})\right).
$$
Here $W_{\xi}$ is the stabilizer of $\xi$ in $W$, while $W_\xi^{\text{aff}}$ is obtained by adding to $W_\xi$ the affine reflection with respect to the hyperplane $(x,\beta_1)=1$, with $\beta_1$ being the largest short root.

Further, let $\mathcal{K}$ be the Kostant cascade of orthogonal roots  for $\Phi$ defined in \cite{K}, which is a collection of positive roots with a partial ordering ``$\preccurlyeq$'' that we recall in Section \ref{KC}. Then we have that

\begin{theorem}\label{Tmain} The Poincar\'e polynomial of $W_\lambda^{\text{aff}}$ factorizes as follows,
\begin{equation}\label{MF}
P\left(W_\lambda^{\text{aff}}/W_\lambda\right) = 1+\sum_{\xi_\beta}t^{\ell_{\xi_\beta}} P\left(W_\lambda/W_{\lambda,\xi_\beta}\right),
\end{equation}
where the summation is over all dominant
\begin{equation}
\xi_\beta = \sum_{\beta' \preccurlyeq \beta} \beta', 
\end{equation} 
with $\beta\in\mathcal{K}$ such that $\langle\xi_\beta,\xi_\beta\rangle=2\langle\xi_\beta,\lambda\rangle$.

\end{theorem}


The formula in Theorem \ref{Tmain} has its origin in the study of the quantization of the moduli stack of parabolic Higgs bundles associated with $\mathbb{P}^1$ with two marked points, which we briefly explain now. 

Let $M_H$ be the moduli stack of $G$-Higgs bundles on a Riemann surface $\Sigma$, and $L$ be the determinant line bundle over $M_H$. The Verlinde formula for Higgs bundles, motivated by quantum physics in \cite{GP,GPYY} and proved in \cite{AGP} and \cite{HL}, gives the graded dimension of $H^0(M_H, L^k)$
$$ \dim_t H^0(M_H, L^k) = \sum_{n=0}^\infty \dim H_n^0(M_H, L^k) t^n.$$
Here, we have used the lift of the ${\mathbb C}^*$-action on $M_H$ to the space of holomorphic sections on $L$ to decompose $H^0(M_H, L^k)$, with $H_n^0(M_H, L^k)$ being a subspace where ${\mathbb C}^*$ acts by the $n$-th tensor power of the defining action of ${\mathbb C}^*$ on ${\mathbb C}$. 

The Verlinde formula for Higgs bundles extends to the case of parabolic Higgs bundles and naturally defines a one parameter (denoted $t$) family of 2d TQFT as  established in  \cite{AGP}. Such a one parameter family of TQFT is uniquely determined by a one parameter family of commutative Frobenius algebras, whose underlying family of finite dimensional vector space  $V_t= V$ is independent of $t$ and has a basis given by the integrable weights of $G$ at level $k$. The associated one parameter family of bilinear forms
$$ B^{(t)}: V^{\otimes 2} \ra \mathbb{C}$$
is determined by applying the TQFT functor to $\mathbb{P}^1$ with two marked points.
Given two integrable weights $\lambda_1$ and $\lambda_2$, one can consider the matrix element $B^{(t)}_{\lambda_1\lambda_2}$, and Theorem~\ref{Tmain} in the present paper was used in \cite{AGP} to prove the following elegant formula
$$
B^{(t)}_{\lambda_1\lambda_2}=\delta_{\lambda_1\lambda_2^*}P\left(W_{\lambda_1}^{\text{aff}}\right).
$$
More precisely, $B^{(t)}_{\lambda_1,\lambda_2}$ is given by an index over the moduli stack of $G$-bundles on $\mathbb{P}^1$ and can be computed by summing over the Shatz strata, labeled by elements in the co-root lattice of $G$ (see Section 5 of \cite{AGP} for details). The main theorem of this paper ensures that contributions from different Shatz strata can be nicely summed together.\footnote{Notice that $\Phi$ in this paper is in fact the coroot lattice of $G$, which results in our convention for $W_\lam^{\text{aff}}$, the inner product $(\cdot,\cdot)$ and $\mathcal{K}$.} 

Furthermore, this interpretation of the main result of this paper also hints at several possible ways to generalize our main theorem. For example, as the Poincar\'e polynomial $P\left(W_{\lambda}\right)$ arises as the normalization factor for the Hall-Littlewood polynomials, it is natural to consider replacing $P\left(W_{\lambda}\right)$ with the normalization factor for other types of orthogonal polynomials associated with root systems. 

Our factorization formula also has the following geometric interpretation. Denote the maximal compact subgroup of $G$ as $K$, whose maximal torus is denoted as $T$ with Cartan subalgebra $\frak{t}$. Given an element $\lambda$ in the Weyl alcove in $\frak{t}$, we can define two subgroups of $K$. The first, $K_\lambda$, is the stabilizer of $\lambda$, while the second, $K'_\lambda$, is the stabilizer of $\exp(\lambda)\in K$. Then their Weyl groups will be $W_\lambda$ and $W_\lambda^{\text{aff}}$ respectively. On the other hand, using the following fibration of (generalized) flag varieties,
$$
K_\lambda/T \rightarrow K'_\lambda/T\rightarrow K'_\lambda/K_\lambda,
$$
we can relate their cohomology groups using the Serre spectral sequence. As all their cohomology groups are in even degrees, all differentials in the spectral sequence vanish, and we have
$$
H^*(K'_\lambda/T)=H^*(K_\lambda/T)\otimes H^*(K'_\lambda/K_\lambda).
$$
This gives a factorization of Poincar\'e polynomials of these flag varieties. As
$$
P(K_\lambda/T)=P(W_\lambda) 
$$
and
$$
P(K'_\lambda/T)=P\left(W_\lambda^{\text{aff}}\right),
$$
we get the following formula for the Poincar\'e polynomial of the partial flag variety $K'_\lambda/K_\lambda$,
$$ 
P(K'_\lambda/K_\lambda)=P\left(W_\lambda^{\text{aff}}\right)/P(W_\lambda).
$$
Therefore, the main theorem of the present paper gives an explicit expression for this Poincar\'e polynomial.

\bigskip


The proof of the main theorem stretches over Section \ref{Gaussian} to Section \ref{proofMT}. A very preliminary remark about Theorem \ref{Tmain} is that if $\lambda$ is not on the affine wall of $\frak{A}$, then (\ref{MF}) is trivially satisfied, since the left-hand side is $1$ and so is the right-hand side, where the sum is over the empty set. 

In Section \ref{Gaussian} we present proofs of elementary identities between Gaussian polynomials, which we use in the following Section \ref{Poincare} to obtain equation (\ref{2.17}), which is a preliminary version of formula (\ref{MF}) from Theorem \ref{Tmain} above in the special case where $\lambda$ is the minuscule fundamental weight associated to a simple root, and the range of the sum on the left-hand side of (\ref{MF}) is not yet understood to be given as stated in Theorem \ref{Tmain}, but is simply given in {\tt Table~2}  for the different types of root systems and the exponents of $t$ are not yet determined.

In Section \ref{StabW} we establish that formula (\ref{2.17}) can be rewritten to formula (\ref{(*)}), which refers to the affine stabilisers as opposed to the full Weyl group in case of (\ref{2.17}), thus bringing it one step closer to (\ref{MF}). A remark following the proof of Lemma \ref{L3.3} establishes that the right-hand side of (\ref{MF}) can be identified for any $\lambda$ on the affine wall with the corresponding right-hand side of (\ref{2.17}) for some minuscule fundamental weight associated to a simple root of some other root system. There remains now to prove that the exponents of $t$ is given as in (\ref{MF}) or equivalently as in (\ref{3.8}) and of course also that the summation range and each coefficient is correct. As is stated in Lemma \ref{L3.4}, by inspecting through {\tt Table 3}, one sees that the exponents which are involved in (\ref{2.17}) are indeed given by (\ref{3.8}) for the case of minuscule fundamental weights. That the coefficient from (\ref{MF}) works is the content of Proposition \ref{P3.6}. But it still leaves the determination of the summation range.

In order to determine a suitable summation range, Section \ref{KC} recalls fundamentals about the Kostant cascade of an indecomposable root system and the following section introduces the set (\ref{5.2}) of dominant weights which constitute the summation range in (\ref{MF}). Proposition \ref{PA.1} and the following five lemmas establishes various properties which are needed in the next section. 

Section \ref{proofMT} begins with Lemma \ref{L5.2}, which establishes that for the minuscule weights associated to simple roots, the summation range given in {\tt Table~2} is indeed given by the set (\ref{5.2}). The following two lemmas provide further properties of the set (\ref{5.2}), which allow us to prove the main theorem in its reformulated version of Theorem \ref{Tmainr} in Section \ref{proofMT}. This is now done by inspection, where we use the result of Section \ref{StabW} to map the problem for general $\lambda$ to a problem concerning the minuscule fundamental weight associated to a simple root of some other root system and then identify the sum on the right hand of (\ref{MF}) with the corresponding sum for this minuscule fundamental weight. This is done by matching up the summation range, the exponents of $t$ and their corresponding coefficients.

\section{Gaussian polynomials}\label{Gaussian}

\bigskip
We work with polynomials over~$\Zet$ in one variable~$t$. 
For any integers $n \geq r \geq 0$ consider the Gaussian polynomial
\begin{equation}
\left[ {n \atop r} \right] := {(t^n - 1) \, (t^{n-1} - 1) \, \ldots \, (t^{n-r+1} - 1)
\over (t^r - 1) \, (t^{r-1} - 1) \, \ldots \, (t -1)}.$$ One has $\left[ {n \atop 0} \right] 
= 1 = \left[ {n \atop n} \right]$ and for $n > 0$ $$\left[ {n \atop r} \right] = 
\left[ {n-1 \atop r} \right] + t^{n-r} \, \left[ {n-1 \atop r-1} \right] \label{1.1}
\end{equation}
where the first summand on the right is to be interpreted as~$0$ if $r = n$
and the second one as~$0$ if~$r=0$.

A classical result (see e.g. Example 3 on page 26 in \cite{McD})
says: If $u$ is another variable, then 
\begin{equation}
\prod_{l=0}^{n-1} 
(1 + t^l \, u) = \sum_{r=0}^n \, t^{r (r-1)/ 2} \, \left[ {n \atop r} \right] \, u^r. 
\label{1.2}
\end{equation} 
This can be proved by induction on~$n$ by using~(\ref{1.1}).

If we plug in $t \, u$ for~$u$ in (\ref{1.2}) and use $r + \frac 1 2 \, r \, (r-1) 
= \frac 1 2 \, r \, (r+1)$, then we get   
\begin{equation}
\prod_{l=1}^{n} (1 + t^l \, u) = 
\sum_{r=0}^n \, t^{r (r+1)/ 2} \, \left[ {n \atop r} \right] \, u^r. \label{1.3}
\end{equation}

\begin{claim} \label{C1}
 If $0 \leq m \leq n$, then
\begin{equation}
\left[ {m + n \atop m} \right] = \sum_{s=0}^m \, t^{s^2} \, 
\left[ {m \atop s} \right] \, \left[ {n \atop s} \right]. \label{1.4}
\end{equation}
\end{claim}

\Proof: We start with~(\ref{1.2}) for~$m+n$: 
\begin{eqnarray*}
\sum_{r=0}^{m+n} \, 
t^{r (r-1)/ 2} \, \left[ {m+n \atop r} \right] \, u^r & = &\prod_{l=0}^{m-1} 
(1 + t^l \, u) \, \prod_{l=m}^{m+n-1} (1 + t^l \, u) \cr & = &\prod_{l=0}^{m-1} 
(1 + t^l \, u) \, \prod_{l=0}^{n-1} (1 + t^l \, t^m u)\cr &= &\left( \sum_{r=0}^m
 t^{r (r-1)/ 2} \left[ {m \atop r} \right]  u^r \right) \left( \sum_{s=0}^n
 t^{s (s-1)/ 2} \left[ {n \atop s} \right] t^{ms}  u^s \right). 
 \end{eqnarray*}
A comparison of the coefficient of~$u^p$ on both sides yields
$$
t^{p (p-1)/2} \, \left[ {m+n \atop p} \right] = \sum_{r+s=p} \, 
t^{{r \choose 2} + {s \choose 2} + m \, s} \,  \left[ {m \atop r} \right]  
 \left[ {n \atop s} \right],
 $$  
 in particular for $p = m$ (using $m \leq n$)
 $$ \left[ {m+n \atop m} \right] = \sum_{s=0}^m \, t^{{m-s \choose 2} + 
 {s \choose 2} + m \, s - {m \choose 2}} \,  \left[ {m \atop m-s} \right]  
 \left[ {n \atop s} \right].$$ A simple calculation shows that 
 ${m-s \choose 2} +  {s \choose 2} + m \, s - {m \choose 2} = s^2$. 
 Now the claim follows from $\left[ {m \atop m-s} \right] = 
 \left[ {m \atop s} \right]$.  

\begin{claim} \label{C2} Let $n$ be a positive integer. Set 
$$
S_n^{\rm ev} = \sum_{0 \leq r \leq n, \, r \, even} \, t^{r (r-1)/ 2} \, 
\left[ {n \atop r} \right], \qquad  S_n^{\rm odd} = \sum_{1 \leq r \leq n, 
\, r \, odd} \, t^{r (r-1)/ 2} \, \left[ {n \atop r} \right].
$$ 
Then
\begin{equation}
S_n^{\rm ev} = S_n^{\rm odd} = \prod_{l=1}^{n-1} (1 + t^l).\label{1.5}
\end{equation}
\end{claim}

\Proof: Note that (\ref{1.2}) evaluated at~$u=1$ yields $$S_n^{\rm ev} + 
S_n^{\rm odd} = 2 \, \prod_{l=1}^{n-1} (1 + t^l),$$ where the~$2$ arises
from the factor $1 + u$ for~$l=0$ in~(\ref{1.2}). Similarly we get at $u=-1$ in~(\ref{1.2}) that
$$S_n^{\rm ev} - 
S_n^{\rm odd} = \prod_{l=0}^{n-1} (1 - t^l) = 0.$$
The claim follows immediately from this.

\section{Poincar\'e polynomials}\label{Poincare}

\bigskip
Let $\Phi$ be a finite root system.
We denote by~$\Delta$ a basis for~$\Phi$ and by~$W$ its Weyl group. 
We shall below use the same numbering as in the tables in \cite{B}
for simple roots~($\alp_i$), simple reflections ($s_i = s_{\alp_i}$) and 
fundamental weights~($\varpi_i)$ of a finite indecomposable root system. 

The Poincar\'e polynomial of a finite Coxeter group such as~$W$ 
will be denoted by $P (W) = \sum_{w \in W} t^{l (w)}$, and we write 
$P (X_n) = P (W)$ if the root system has type~$X_n$.  The degree
of~$P (W)$ is the maximal length in~$W$, equal to the number of 
positive roots in~$\Phi$. One has in general 
\begin{equation}
P (W) = 
\prod_{i=1}^n \, {t^{m_i + 1} - 1 \over t - 1} \label{2.1}
\end{equation}
where $m_1, 
m_2, \ldots, m_n$ are the exponents of the root system. One gets 
explicitly for the classical types 
\begin{equation}
P (A_n) = \prod_{i=1}^n {t^{i+1} - 1 
\over t-1},  \qquad n \geq 1\label{2.2}
\end{equation}
\begin{equation}
P (B_n) = P (C_n) = 
\prod_{i=1}^n {t^{2i} - 1 \over t-1}, \qquad n \geq 1\label{2.3}
\end{equation}
and 
\begin{equation}
P (D_n) = {t^n - 1 \over t -1} \,\prod_{i=1}^{n-1} \, {t^{2i} - 1 \over t-1},
\qquad n \geq 2\label{2.4}
\end{equation} 
where we use the convention that
$B_1 = C_1 = A_1$ and $D_3 = A_3$ and $D_2 = A_1 \times A_1$.
It will be convenient to set $P (A_0) = P (C_0) = 1$ and even 
$P (A_{-1}) = 1$.

We shall need also 
\begin{equation}
P (E_6) = {(t^2 - 1) \, (t^5 - 1) \, (t^6 - 1) \, (t^8 - 1) \, (t^9 -1) 
\, (t^{12} - 1) \over (t-1)^6} \label{2.5}
\end{equation}
and
\begin{equation}
P (E_7) = {(t^2 - 1) \, (t^6 - 1) \, (t^8 - 1) \, (t^{10} - 1) \, (t^{12} -1) 
\, (t^{14} - 1) \, (t^{18} - 1) \over (t-1)^7}. \label{2.6}
\end{equation}

Note: If the root system is a direct sum of several root systems, 
then its Poincar\'e polynomial is the product of the 
Poincar\'e polynomials of its summands.

If $I$ is a subset of~$\Delta$, then we denote by $W_I$ the 
subgroup of~$W$ generated by all $s_i$ with $\alp_i \in I$.
The set $\Phi_I = \Phi \cap \Zet \, I$ is a root subsystem of~$\Phi$
with basis~$I$ and Weyl group~$W_I$. We set $P (W/W_I) = 
P (W) / P (W_I)$. This quotient is in fact  a polynomial, equal to 
$\sum_w t^{\ell (w)}$ where we sum over minimal length
representatives for the cosets of~$W_I$ in~$W$. The degree 
of $P (W / W_I)$ is equal to the number of positive roots in~$\Phi$ 
not in~$\Phi_I$. 

Assume now that $\Phi$ is indecomposable.
We are now going to exhibit some explicit expressions for $P (W/W_I)$ 
in case $I = \Delta \setminus \{ \alp \}$ where the simple root~$\alp$
has the property that the corresponding fundamental weight is 
minuscule. 

We start in type~$A_n$. Here any fundamental weight is minuscule. 
If  $\alp = \alp_i$ with $1 \leq i \leq n$, then the remaining simple roots 
span a root system of type $A_{i-1} \times A_{n-i}$. We get then (as is 
well-known) 
\begin{eqnarray}
{P (A_n) \over P (A_{i-1}) \cdot 
P (A_{n-i})} & = & {\prod_{s=2}^{n+1} (t^s - 1) \over (t-1) \, 
\prod_{p=2}^{i} (t^p - 1) \, \prod_{q=2}^{n-i+1} (t^q - 1)}\\ & = &
{\prod_{s=n+2-i}^{n+1} (t^s - 1) \over \prod_{p=1}^{i} (t^p - 1)} \cr& 
= &\left[ {n+1 \atop i} \right]. \label{2.7}
\end{eqnarray}
Note that this also works for $i = 1$ and $i = 0$ thanks to our 
convention. For $i = 1$ we get more explicitly 
$$
{P (A_n) \over 
P (A_{n-1})} = {t^{n+1} - 1 \over t - 1}.
$$ 
Since 
$$
{t^{n+1} - 1 \over t - 1} 
= 1 + t \, {t^n - 1 \over t - 1}
$$ 
we get thus 
\begin{equation}
{P (A_n) \over 
P (A_{n-1})} = 1 + t \, {P (A_{n-1}) \over P (A_{n-2})}, \qquad 
n \geq 1. \label{2.8}
\end{equation} 
On the other hand, we get from Claim~\ref{C1} 
for $0 \leq m \leq n$ 
\begin{equation}
{P (A_{m+n-1}) \over P (A_{m-1}) \cdot 
P (A_{n-1})} = \sum_{s=0}^m \, t^{s^2} \, {P (A_{m-1}) \cdot P (A_{n-1}) 
\over P (A_{s-1}) \, P (A_{m-s-1}) \, P (A_{s-1}) \, P (A_{n-s-1})}.
\label{2.9}
\end{equation} 
(For $m = 1$ we get (\ref{2.8}) back.)

\medskip
In type~$B_n$ our only choice is $\alp= \alp_n$. Then the 
remaining simple roots span a root system of type $A_{n-1}$. 
We get then 
\begin{equation}
{P (B_n) \over P (A_{n-1})} = {\prod_{i=1}^n 
(t^{2i} - 1) \over (t-1) \, \prod_{i=2}^n (t^i - 1)} = \prod_{i=1}^n
(t^i + 1).$$ Now (\ref{2.7}) and (\ref{1.3}) yield $${P (B_n) \over P (A_{n-1})} 
= \sum_{s=0}^n \, t^{s (s+1)/2} \, {P (A_{n-1}) \over P (A_{s-1}) 
\cdot P (A_{n-s-1})}, \qquad n \geq 2. \label{2.10}
\end{equation}

In type~$C_n$ our only choice is $\alp= \alp_1$. Then  the 
remaining simple roots span a root system of type $C_{n-1}$. 
We get then 
$$
{P (C_n) \over P (C_{n-1})} = {\prod_{i=1}^n 
(t^{2i} - 1) \over (t-1) \, \prod_{i=1}^{n-1} (t^{2i} - 1)} = {t^{2n} - 1
\over  t-1}, \qquad n \geq 1.
$$ 
Since 
$$
{t^{2n} - 1 \over t - 1} = 
1 + t \, {t^{2 (n-1)} - 1 \over t - 1} + t^{2n-1}
$$ 
we get thus 
\begin{equation}
{P (C_n) \over P (C_{n-1})} = 1 + t \, {P (C_{n-1}) \over P (C_{n-2})} 
+  t^{2n-1}, \qquad n \geq 2. \label{2.11}
\end{equation}

\medskip
In type~$D_n$ there are three choices for~$\alp$. We take first $\alp
= \alp_1$; here the remaining simple roots span a root system of 
type~$D_{n-1}$. We get then 
\begin{eqnarray*}
{P (D_n) \over P (D_{n-1})} & 
= {(t^n - 1) \, \prod_{i=1}^{n-1} (t^{2i} - 1) \over (t-1) \, (t^{n-1} - 1) \, 
\prod_{i=1}^{n-2} (t^{2i} - 1)} = {(t^n-1) \, (t^{2(n-1)} - 1)
\over ( t-1) \, (t^{n-1} - 1)} \cr & = {(t^n-1) \, (t^{n-1} + 1)
\over t-1}, \qquad n \geq 3.
\end{eqnarray*} 
Since 
$$
{(t^n-1) \, (t^{n-1} + 1) \over t-1}  
= 1 + t \, {(t^{n-1}-1) \, (t^{n-2} + 1) \over t-1} + t^{2(n-1)}
$$ 
we get 
thus 
\begin{equation}
{P (D_n) \over P (D_{n-1})} = 1 + t \, {P (D_{n-1}) \over 
P (D_{n-2})} +  t^{2(n-1)}, \qquad n \geq 4. \label{2.12}
\end{equation}

If we take $\alp = \alp_{n-1}$ or $\alp = \alp_n$ in type~$D_n$,  then 
the remaining simple roots span a root system of type $A_{n-1}$. 
We get then 
$$
{P (D_n) \over P (A_{n-1})}  = {(t^n - 1) \, 
\prod_{i=1}^{n-1} (t^{2i} - 1) \over (t-1) \, \prod_{i=2}^{n} (t^{i} - 1)} 
= \prod_{i=1}^{n-1} {t^{2i} - 1 \over t^i-1} =\prod_{i=1}^{n-1}
(t^i + 1).
$$ 
Therefore Claim~\ref{C2} says 
\begin{multline} {P (D_n) \over P (A_{n-1})}  
 = \sum_{0 \leq r \leq n, \, r \, \text{even}} t^{r (r-1) / 2} {P (A_{n-1}) \over 
P (A_{r-1}) \cdot P (A_{n-r-1})} \\ = \sum_{0 \leq 2s \leq n} t^{s (2s-1)} 
{P (A_{n-1}) \over P (A_{2s-1}) \cdot P (A_{n-2s-1})}, \qquad n \geq 4.
\label{2.13}
\end{multline}

\medskip
In type $E_6$ we can take $\alp = \alp_1$ or $\alp = \alp_6$. In both 
cases the remaining simple roots span a root system of type $D_5$. 
In type $E_7$ only $\alp = \alp_7$ is possible; in this case the 
remaining simple roots span a root system of type $E_6$. 
The corresponding quotients are 
$$
{P (E_6) \over P (D_5)} 
= {(t^9 - 1) \, (t^{12} - 1) \over (t - 1) \, (t^4 - 1)} = 
{(t^9 - 1) \, (t^8 + t^4 + 1) \over (t - 1)}
$$ 
and 
$${P (E_7) \over P (E_6)} 
= {(t^{10} - 1) \, (t^{14} - 1) \, (t^{18} - 1) \over (t - 1) \, (t^5 - 1) \, (t^9 - 1)} = 
{(t^5 + 1) \, (t^{14} - 1) \, (t^9 + 1) \over (t - 1)}.
$$ 
A little calculation 
shows now that  
\begin{equation}
{P (E_6) \over P (D_5)} = 1 + t \, {P (D_5) \over P (A_4)} 
+ t^8 \, {P (D_5) \over P (D_4)} \label{2.14}
\end{equation}
and  
\begin{equation}
{P (E_7) 
\over P (E_6)} = 1 + (t + t^{10}) \, {P (E_6) \over P (D_5)} + t^{27}.
\label{2.15}
\end{equation}

\bigskip
Here is a table listing the possible pairs $(\Phi, \alp)$ together 
with the type of $\Delta \setminus \{ \alp \}$ and the degree 
of $P (W/W_I)$. 

\begin{center}
\begin{tabular}{ c c c c c }
$\Phi$ & $\alp$ & type $\Delta \setminus \{ \alp \}$ &degree \\\hline
$A_n$& $\alp_j$ & $A_{j-1} \times A_{n-j}$ & $j \, (n+1-j)$ \\
$B_n$ & $\alp_n$ & $A_{n-1}$ & $n \, (n+1)/2$ \\
$C_n$ & $\alp_1$ & $C_{n-1}$ & $ 2\, n-1$ \\
$D_n$ & $\alp_1$ & $D_{n-1}$ & $2 \, (n-1)$ \\
$D_n$ & $\alp_{n-1}$, $\alp_n$ & $A_{n-1}$ & $n \, (n-1)/2$ \\
$E_6$ & $\alp_1$, $\alp_6$ & $D_5$ & $16$ \\
$E_7$ & $\alp_7$ & $E_6$ & $27$
\end{tabular}
$$\hbox{\tt Table 1}$$
\end{center}

We can give the equations (\ref{2.8}) -- (\ref{2.15}) a uniform appearance.
For any dominant (not necessarily integral) element~$\lam$ in the real 
span~$\Err \Phi$ of the roots set 
\begin{equation}
\Delta_\lam = \{\, \alp \in \Delta \mid 
s_\alp \, \lam = \lam \,\} = \{\, \alp \in \Delta \mid (\lam, \alp^\vee) = 0 \,\}.
\label{2.16}
\end{equation}  
Then the stabiliser~$W_\lam$ of~$\lam$ in~$W$ is equal 
to~$W_{\Delta_\lam}$. We have, for example, $\Delta_{\varpi_\alp} 
= \Delta \setminus \{ \alp \}$ for any $\alp \in \Delta$. 
So the denominators on the left hand sides of our equations are just $P(W_{\varpi_\alp})$.

Note that the summands for $s = 0$ in (\ref{2.9}), (\ref{2.10}), and~(\ref{2.13}) are 
equal to~$1$. The denominators on the right hand side of our equations 
can be interpreted as Poincar\'e polynomials of intersections
$W_{\varpi_\alp, \xi} = W_{\varpi_\alp} \cap W_\xi$ 
for suitable~$\xi$. We get thus with the above notations 
\begin{equation}
P (W / W_{\varpi_\alp}) = 1 + \sum_{\xi \in \CM} \, t^{\ell_\xi} \, 
P (W_{\varpi_\alp} / W_{\varpi_\alp, \xi}) \label{2.17}
\end{equation}   
where the $\ell_\xi$ 
are suitable positive integers and the set~$\CM$ is given by {\tt Table~2}
below (where we skip two cases that are symmetric to some we include).
\vskip.3cm
\begin{center}
\begin{tabular}{ c c c c c }
$\Phi$ & $\alp$ & $\CM$\\ \hline
$A_n$& $\alp_j$  & $\varpi_s + \varpi_{n+1-s}$, $1 \leq s \leq \min (j, n-j+1)$\\
$B_n$ & $\alp_n$ & $\varpi_s$, $1 \leq s < n$, $2 \varpi_n$ \\
$C_n$ & $\alp_1$ & $\varpi_2$, $2 \varpi_1$\\
$D_n$ & $\alp_1$ & $\varpi_2$, $2 \varpi_1$\\
$D_n$ & $\alp_n$ & if $n$ even: $\varpi_{2s}$, $1 \leq s < n/2, \, 2 \varpi_n$\\
&&if $n$ odd: $\varpi_{2s}$, $1 \leq s < (n-1)/2, \, \varpi_{n-1} + \varpi_n$\\
$E_6$ & $\alp_6$ & $\varpi_2$, $\varpi_1 + \varpi_6$\\
$E_7$ & $\alp_7$ &$\varpi_1$, $\varpi_6$, $2 \varpi_7$\\
\end{tabular}
\end{center}
$$\hbox{\tt Table 2}$$

Note that these elements are not uniquely determined by~(\ref{2.17}).
We shall see that there is a general recipe that produces exactly 
the set~$\CM$ above and also find a general expression 
for the exponents~$\ell_\xi$.

\section{Stabilisers in the affine Weyl group}
\label{StabW}

\bigskip

Assume from now on that $\Phi$ is indecomposable.
We consider here the affine Weyl group~$W^{\rm aff}$ of the root 
system as the group generated by~$W$ and all translations by roots
(not by coroots). Denote the largest short root in~$\Phi$ by~$\beta_1$. 
(Our convention is that all roots are short if all roots have the same length.) 
Set $s_0 = s_{\beta_1, 1} \in W^{\rm aff}$ equal to the affine reflection 
with respect to the hyperplane $(x, \beta_1^\vee) = 1$. Now
$W^{\rm aff}$ is a Coxeter group with Coxeter generators 
$$S^{\rm aff} = \{\, s_\alp \mid \alp \in \Delta \,\}  \cup \{ s_0 \}.$$ 

Set $\alp_0 = - \beta_1$. The extended Dynkin diagram of~$\Phi$ has
vertices corresponding to $\Delta \cup \{ \alp_0\}$ and is 
constructed after the usual rules. In particular, the vertex corresponding
to~$\alp_0$ is linked to some $\alp \in \Delta$ if and only if 
$(\alp_0, \alp) < 0$ --- equivalently: if and only if $(\beta_1, \alp) > 0$ 
--- and the type of the link is determined by the fact that $\alp_0$
is short. The vertices in $\Delta \cup \{ \alp_0\}$ are in bijection with
the elements in~$S^{\rm aff}$, and the Coxeter graph of the
Coxeter system~$(W^{\rm aff}, S^{\rm aff})$ is the Coxeter graph 
associated to the extended Dynkin diagram.

The fundamental alcove 
\begin{equation}
\frak A = \{\, \mu \in \Err \Phi \mid 
\hbox {$0 \leq (\mu, \alp_i^\vee)$ 
for all $i$, $1 \leq i \leq n$, and 
$(\mu, \beta_1^\vee) \leq 1$} \,\} \label{3.1}
\end{equation}
is a fundamental domain 
for~$W^{\rm aff}$. For any $\mu \in {\frak A}$ the stabiliser%
~$W_\mu^{\rm aff}$ of~$\mu$ in~$W^{\rm aff}$ is generated by all 
$s \in S^{\rm aff}$ with $s \mu = \mu$.

Any proper subset $K$ of $\Delta \cup \{ \alp_0 \}$ is the Dynkin 
diagram of a finite root system. We denote its Weyl group by~$W_K$
and identify it with the subgroup of~$W^{\rm aff}$ generated by 
all $s_\alp$ with $\alp \in K \cap \Delta$ together with~$s_0$ in 
case $\alp_0 \in K$. (For $K \subset \Delta$ this is compatible 
with our earlier definition.) 

The statement above on stabilisers in~$W^{\rm aff}$ for elements 
in~$\Ahh$ implies with the notation from~(\ref{2.16}) 
\begin{equation}
\hbox {$\lam \in 
\Ahh$ and $(\lam, \beta_1^\vee) = 1$} \; \Longrightarrow \; 
W^{\rm aff}_\lam = W_{\Delta_\lam \cup \{ \alp_0\}}. \label{3.2}
\end{equation}

\begin{lemma} \label{L3.1} Let $\alp \in \Delta$ be a simple root such that 
$\varpi_\alp$ is a minuscule fundamental weight. Then there exists 
an automorphism of the extended Dynkin diagram that interchanges 
$\alp$ and~$\alp_0$. If we remove $\alp$ from the extended diagram, 
then the remainder $(\Delta \setminus \{ \alp \}) \cup \{ \alp_0 \}$ 
is isomorphic to~$\Delta$, and we get the extended diagram of 
the remainder by adding~$\{ \alp \}$. 
\end{lemma}

\Proof: By inspection and well-known. 

\medskip
Recall that a fundamental weight~$\varpi_\alp$ is minuscule if and only 
if $1 = (\varpi_\alp, \beta_1^\vee)$. So we have $\varpi_\alp 
\in \Ahh$ and can apply (\ref{3.2}) to it. 

\begin{corollary} \label{Co3.2} Let $\alp \in \Delta$ be a simple root such that 
$\varpi_\alp$ is a minuscule fundamental weight. Then $W^{\rm aff}_{
\varpi_\alp}$ is isomorphic to~$W$ as a Coxeter group. 
\end{corollary}

\Proof: Since $\Delta_{\varpi_\alp} = \Delta \setminus \{ \alp \}$ we get 
from~(3.2) that $W^{\rm aff}_{\varpi_\alp} = W_{(\Delta \setminus \{ \alp \})
\cup \{ \alp_0 \}}$. An automorphism as in the lemma maps 
$(\Delta \setminus \{ \alp \}) \cup \{ \alp_0 \}$ onto~$\Delta$. Therefore 
$W^{\rm aff}_{\varpi_\alp}$ is isomorphic to~$W_\Delta = W$ as a 
Coxeter group. 

\Remark: We get here in particular that $P (W^{\rm aff}_{\varpi_\alp}) 
= P (W)$. So we can rewrite (\ref{2.17}) as 
\begin{equation}
P (W^{\rm aff}_{\varpi_\alp} 
/ W_{\varpi_\alp}) = 1 + \sum_{\xi \in \CM} \, t^{\ell_\xi} \, 
P (W_{\varpi_\alp} / W_{\varpi_\alp, \xi}). \label{(*)}
\end{equation} 
Our goal 
is to generalise this formula so that we can replace~$\varpi_\alp$ 
by any~$\lam \in \Ahh$ satisfying $(\lam, \beta_1^\vee) = 1$. 

\medskip
For any dominant~$\lam$ as in~\eqref{3.1} denote by $J (\lam) \subset 
\Delta_\lam$ the subset such that 
\begin{equation}
J (\lam) \cup \{ \alp_0 \} = 
\hbox {the connected component of~$\Delta_\lam \cup \{\alp_0 \}$ 
containing~$\alp_0$}. \label{3.3}
\end{equation}
Set $I (\lam) = \Delta_\lam 
\setminus J (\lam)$. 

For example, we have $J (\varpi_\alp) = \Delta \setminus \{ \alp \}$ 
and $I (\varpi_\alp) = \emptyset$ for any $\alp \in \Delta$ with 
$\varpi_\alp$ minuscule as 
$( \Delta \setminus \{ \alp \}) \cup \{ \alp_0 \}$ is connected 
(isomorphic to~$\Delta$). On the other hand, we get $J (\lam) = 
\emptyset$ if and only if $\alp \notin \Delta_\lam$ for all $\alp \in \Delta$ 
with $(\beta_1, \alp^\vee) > 0$.  This means in particular $J (\beta_1) 
= \emptyset$ and $I (\beta_1) = \Delta_{\beta_1}$.

For any dominant~$\lam$ as above we get a direct product decomposition 
\begin{equation}
W_\lam = W_{J (\lam)} \times  W_{I (\lam)} \label{3.4}
\end{equation}
since $\Delta_\lam$ is the disjoint union of $J (\lam)$ and $I (\lam)$, 
and since no root in $J (\lam)$ is linked to a root in~$I (\lam)$ by 
definition of~$J (\lam)$.

If $\lam \neq 0$, then $\Delta_\lam \cup \{ \alp_0 \}$ and $J (\lam) 
\cup \{ \alp_0 \}$ are proper subsets of the extended Dynkin diagram 
and define finite parabolic subgroups $W_{\Delta_\lam \cup \{ \alp_0 \}}$
and $W_{J (\lam) \cup \{ \alp_0 \}}$ of~$W^{\rm aff}$. We get 
similarly to~(\ref{3.4}) that 
\begin{equation}
W_{\Delta_\lam \cup \{ \alp_0 \}} = 
W_{J (\lam) \cup \{ \alp_0 \}} \times  W_{I (\lam)}.\label{3.5}
\end{equation}  
A comparison of (\ref{3.4}) and (\ref{3.5}) implies for the Poincar\'e polynomials 
\begin{equation}
P (W_{\Delta_\lam \cup \{ \alp_0 \}} / W_\lam) = P (W_{J (\lam) \cup 
\{ \alp_0 \}} / W_{J (\lam)}). \label{3.6}
\end{equation}
The following lemma applies to any set of the form
$J (\lam) \cup \{ \alp_0 \}$ where $\lam \neq 0$.

\begin{lemma} \label{L3.3} Let $K$ be a connected proper subset of the 
extended Dynkin diagram $\Delta \cup \{ \alp_0 \}$ such that 
$\alp_0 \in K$. Then $\alp_0$ corresponds to a minuscule 
fundamental weight of the root system with Dynkin diagram~$K$. 
\end{lemma}

\Proof: The claim is obvious if $K$ has type~$A_r$ for some~$r$ 
since here all fundamental weights are minuscule. This takes 
care of $\Phi$ of type~$A_n$ where any possible~$K$ has type~$A$.

So assume that $\Phi$ is not of $A$--type. Then $\alp_0$ is an 
end-vertex of the extended Dynkin diagram $\Delta \cup \{ \alp_0 \}$,
hence also of~$K$. Furthermore $\alp_0$ corresponds to a short 
root. Now in the classical cases (BCD-types) any short simple root 
located at an end of the Dynkin diagram corresponds to a minuscule 
fundamental weight. 

So we are left with the possibility that $K$ has exceptional type. 
If $K$ has type $E_8$, $F_4$, or~$G_2$, then $K$ can occur 
as a subdiagram in $\Delta \cup \{ \alp_0 \}$ only if $\Delta$ has the same
type as~$K$, and it has to be the subdiagram $\Delta$ in $\Delta \cup \{ \alp_0 \}$.
But we assume that $\alp_0 \in K$. So these cases cannot occur. 

Let us now suppose that $K$ has type~$E_6$. Then $\Phi$ has to 
have type~$E_n$ with $n \geq 6$. Then $\alp_0$ has distance
(in an obvious sense) at least~$2$ from the branching point in 
$\Delta \cup \{ \alp_0 \}$, hence also at least distance~$2$ from the 
branching point in~$K$. Therefore $\alp_0$ corresponds to one 
of the two minuscule fundamental weights.

The argument is similar for $K$ of type~$E_7$. Now $\Phi$ has to 
have type~$E_n$ with $n \geq 7$. Then $\alp_0$ has distance
at least~$3$ from the branching point in $\Delta \cup \{ \alp_0 \}$, 
hence also at least distance~$3$ from the branching point in~$K$. 
Therefore the end-vertex $\alp_0$ corresponds to the minuscule 
fundamental weight.

\Remark: Consider $\lam \neq 0$ as above. The lemma says that 
there exists an indecomposable finite root system~$\Phi'$ with 
basis~$\Delta'$ and root $\alp' \in \Delta'$ corresponding to  a 
minuscule fundamental weight such that the pair of Dynkin
diagrams $(J (\lam) \cup \{ \alp_0 \}, J (\lam))$ is isomorphic 
to the pair of Dynkin diagrams $(\Delta', \Delta' \setminus \{ \alp' \})$.
It then follows that $P (W_{J (\lam) \cup \{ \alp_0 \}} / W_{J (\lam)})
= P (W_{\Delta'} / W_{\Delta' \setminus \{ \alp' \}})$ is one of the 
left hand sides considered in~\eqref{2.17}. If we assume in addition 
that $\lam \in \Ahh$ with $(\lam, \beta_1^\vee) = 1$, then we have 
by \eqref{3.3} and~\eqref{3.6} that 
\begin{equation}
P (W^{\rm aff}_\lam / W_\lam) = P (W_{\Delta'} 
/ W_{\Delta' \setminus \{ \alp' \}}).\label{3.7}
\end{equation} 
The question is then 
whether we can also relate the right hand side in~\eqref{2.17} to~$\lam$.

\medskip
With a view to this goal we first want to give a general formula for 
the exponents~$\ell_\xi$ in~(\ref{2.17}). We set for any dominant~$\xi \neq 0$ 
\begin{equation}
\ell_\xi = \deg P (W_{J(\xi) \cup \{ \alp_0 \}} / W_{J (\xi)}), \label{3.8}
\end{equation} 
which obviously agrees with our definition of $\ell_\xi$ in the introduction. Note that we can apply the considerations above also to~$\xi$, 
hence have an isomorphism of pairs of Dynkin diagrams of the
form $(J (\xi) \cup \{ \alp_0 \}, J (\xi)) \simeq (\Delta', \Delta' \setminus 
\{ \alp' \})$ and get then $\ell_\xi = \deg P (W_{\Delta'} 
/ W_{\Delta' \setminus \{ \alp' \}})$; the latter degree can then  be 
read off {\tt Table~1}. The following table gives a list of these 
types and the corresponding~$\ell_\xi$  for certain~$\xi$ that will turn out 
to play a role later on.  

\vskip.3cm
\begin{center}
\begin{tabular}{ c c c c c }
type $\Phi$ & $\xi$ & type $J (\xi)$ & type $J (\xi) \cup \{ \alp_0 \}$ & $\ell_\xi$\\ \hline
any & $\beta_1$ & $\emptyset$ & $A_1$ & $1$\\
$A_n$ & $\varpi_i + \varpi_{n+1-i}, \, 1 < i \leq n/2$ & $A_{i-1} \times A_{i-1}$ & $A_{2i-1}$ & $i^2$\\
$B_n$ & $\varpi_i, \, 1 < i <  n$ & $A_{i-1}$ & $B_i$ & $i (i+1) /2$\\
$B_n$ & $2 \, \varpi_n$ & $A_{n-1}$ & $B_n$ & $n (n+1) /2$\\
$C_n$ & $\varpi_{2i}, \, 1 < i \leq n/2$ & $A_{2i-1}$ & $D_{2i}$ & $i (2i-1)$\\
$C_n$ & $2 \, \varpi_1$ & $C_{n-1}$ & $C_n$ & $2\, n - 1$\\
$D_n$ & $\varpi_{2i}, \, 1 < i < (n-1)/2$ & $A_{2i-1}$ & $D_{2i}$ & $i (2i-1)$\\
$D_n$ & $2 \, \varpi_{n-1}$ or $2 \, \varpi_n$& $A_{n-1}$ & $D_n$ & $n(n-1)/2$\\
$D_{2m+1}$ & $\varpi_{2m} + \varpi_{2m+1}$ & $A_{2m-1}$ & $D_{2m}$ & $m(2m-1)$\\
$D_n$ & $2 \, \varpi_1$ & $D_{n-1}$ & $D_n$ & $2\, (n-1)$\\
$E_6$ & $\varpi_1 + \varpi_6$ & $D_4$ & $D_5$ & $8$\\
$E_7$ & $\varpi_6$ & $D_5$ & $D_6$ & $10$\\
$E_7$ & $2 \, \varpi_7$ & $E_6$ & $E_7$ &$27$\\
$E_8$ & $\varpi_1$ & $D_7$ & $D_8$ & $14$\\
$F_4$ & $\varpi_1$ & $C_3$ & $C_4$ & $7$\\
\end{tabular}
\end{center}
$$\hbox{\tt Table 3}$$

\begin{lemma} \label{L3.4} 
Let $\alp \in \Delta$ be a simple root such that 
$\varpi_\alp$ is a minuscule fundamental weight. Then the 
exponents~$\ell_\xi$ in~$(\ref{2.17})$ coincide with the exponents defined 
by~\eqref{3.8}. We have $\alp \notin J (\xi)$ for all~$\xi \in \CM$.
\end{lemma}

\Proof: Each $\xi$ from {\tt Table~2} occurs also in {\tt Table~3}. One 
can then easily compare the data in~{\tt Table~3} with the exponents 
in (\ref{2.8}) -- (\ref{2.15}). (Note that $\beta_1$ occurs in different guises 
in {\tt Table~2}.)

One checks that $\alp \notin J (\xi)$ by inspection.

\medskip
Let $\lam \in \Ahh$ with $(\lam, \beta_1^\vee) = 1$. In our intended 
generalisation of~\eqref{(*)} to~$\lam$, the right hand side is supposed to 
involve terms of the form $P (W_\lam / W_{\lam, \xi})$. We shall have 
to relate $W_{\lam, \xi}$ to~$J (\lam)$.

\begin{lemma} \label{L3.5} Let $\lam \in \Ahh$ with $(\lam, \beta_1^\vee) = 1$. 
Let $\xi \in \Err \Phi$, $\xi \neq 0$, be dominant. If $I (\lam) \subset 
\Delta_\xi$, then $P (W_\lam / W_{\lam, \xi} ) = P (W_{J (\lam)} / W_{J (\lam) 
\cap \Delta_\xi})$.
\end{lemma}

\Proof: The assumption  $I (\lam) \subset \Delta_\xi$ implies $$\Delta_\lam 
\cap \Delta_\xi = (J (\lam) \cup I (\lam))  \cap \Delta_\xi  = (J (\lam) \cap 
\Delta_\xi) \cup I (\lam).$$ Since $W_{\lam, \xi} =  W_\lam \cap W_\xi$ is 
the parabolic subgroup of~$W$ generated by the $s_\alp$ with $\alp \in 
\Delta_\lam \cap \Delta_\xi$ we get as for~(\ref{3.4}) that $W_{\lam, \xi} = 
W_{ J (\lam) \cap \Delta_\xi} \times W_{I(\lam)}$. A comparison with
(\ref{3.4}) then yields the claim.

\begin{proposition}  \label{P3.6} Let $\lam \in \Ahh$ with $(\lam, \beta_1^\vee) 
= 1$. There exists a finite set~$\CM$ of non-zero dominant weights
such that 
\begin{equation}
P (W^{\rm aff}_\lam / W_\lam) = 1 + \sum_{\xi \in \CM} 
\belowdisplayskip=0pt
\, t^{\ell_\xi} \, P (W_\lam / W_{\lam, \xi}). \label{3.9}
\end{equation} 
\end{proposition}


\Proof:  Because of (\ref{3.2}), (\ref{3.6}), and Lemma~\ref{L3.5} we actually want
to find $\CM$ with 
\begin{equation}
I (\lam) \subset \Delta_\xi \qquad \hbox 
{for all $\xi \in \CM$} \label{3.10}
\end{equation} 
such that 
\begin{equation}
P (W_{J (\lam) \cup \{ \alp_0 \}} / W_{J (\lam)}) = 1 + \sum_{\xi \in \CM} 
\, t^{\ell_\xi} \, P (W_{J (\lam)} / W_{J (\lam) \cap \Delta_\xi}). \label{3.11}
\end{equation} 

Consider $J (\lam) \cup \{ \alp_0 \}$ as the Dynkin diagram of a 
finite root system and take its extended Dynkin diagram. It contains 
an additional vertex~$\alp'$. Set $\Delta' = J (\lam) \cup \{ \alp' \}$
and regard $\Delta'$ as the basis of a  finite root system~$\Phi'$.
By Lemma~\ref{L3.3}, $\alp'$ corresponds to a minuscule fundamental 
weight of~$\Phi'$, so we can apply (\ref{2.17}) to $\Phi'$ and~$\alp'$.
Lemma~\ref{L3.1} implies that we have an isomorphism of pairs of Dynkin 
diagrams 
$$
(J (\lam) \cup \{ \alp_0 \}, J (\lam)) \isto (\Delta', 
\Delta' \setminus \{ \alp'\}).
$$ 
Therefore we can rewrite (\ref{2.17}) as 
\begin{equation}
P (W_{J (\lam) \cup \{ \alp_0 \}} / W_{J (\lam)}) = 1 + \sum_{\xi' \in \CM'} 
\, t^{\ell_{\xi'}} \, P (W_{J (\lam)} / W_{J (\lam) \cap \Delta'_{\xi'}}) \label{3.12}
\end{equation} 
where $\CM'$ is now  a set of dominant weights for~$\Phi'$ and
where $\Delta'_{\xi'}$ is defined as in~(\ref{2.16}), just working in~$\Phi'$ 
instead of~$\Phi$. We now need the following lemma.

\begin{lemma} \label{L3.7}  
There exists a finite set~$\CM$ of non-zero 
dominant weights together with a bijection $\CM \isto \CM'$, 
$\xi \mapsto \xi'$, such that 
\begin{equation}
I (\lam) \subset \Delta_\xi, \qquad 
J (\xi) \subset J (\lam), \qquad \hbox {and} \qquad J (\lam) \cap \Delta_\xi = 
J (\lam) \cap \Delta'_{\xi'} \label{3.13}
\end{equation} 
for all $\xi \in \CM$. 
\end{lemma}


Let us postpone the proof of the lemma and show that (\ref{3.11}) holds 
for any~$\CM$ satisfying this lemma. 

To begin with, the third condition in~(\ref{3.13}) says that the Poincar\'e 
polynomials on the right hand side of~(\ref{3.12}) are the same as
those on the right hand side of~(\ref{3.11}). Thanks to  the first condition 
in~(\ref{3.13}) they are also the same as  in~(\ref{3.9}). 

So it is left to check that $\ell_\xi = \ell_{\xi'}$ for all~$\xi \in \CM$
which means 
\begin{equation}
\deg P (W_{J (\xi) \cup \{ \alp_0 \}} / W_{J (\xi)}) = 
\deg P (W_{J' (\xi') \cup \{ \alp_0 \}} / W_{J' (\xi')}). \label{3.14}
\end{equation} 
Here $J' (\xi') \subset \Delta' = J (\lam) \cup \{ \alp' \}$ is defined 
analogously so that $J' (\xi') \cup \{ \alp_0 \}$ is the connected 
component of $\Delta'_{\xi'} \cup \{ \alp_0 \}$ containing~$\alp_0$.
(Note that $\alp_0$ is by Lemma~\ref{L3.1} also the extra vertex in the 
extended Dynkin diagram of~$\Phi'$.) 

Lemma~\ref{L3.4} says that $\alp' \notin J' (\xi')$, hence $J' (\xi') \subset 
J (\lam)$. So what we really want is that 
\begin{equation}
J (\xi) = J' (\xi') \label{3.15}
\end{equation}
which then implies (\ref{3.14}). Since $$J' (\xi') \subset \Delta'_{\xi'} \cap J (\lam) = 
\Delta_\xi \cap J (\lam) \subset \Delta_\xi$$ thanks to~(\ref{3.13}), the connected 
subset $J' (\xi') \cup \{ \alp_0 \}$ is contained in $\Delta_\xi \cup \{ \alp_0 \}$,
hence in its connected component $J (\xi) \cup \{ \alp_0 \}$. 
It follows that $J' (\xi') \subset J (\xi)$. 

Similarly 
$$J (\xi) \subset \Delta_{\xi} \cap J (\lam) = 
\Delta'_{\xi'} \cap J (\lam) \subset \Delta'_{\xi'}$$ 
implies by the same 
arguments $J (\xi) \subset J' (\xi')$. 

So the proposition follows as soon as we have proved the lemma.


\noindent {\it Proof of Lemma \ref{L3.7}:} Any $\xi \in \Err \Phi$ can be written as $\xi = 
\sum_{\gamma \in \Delta} (\xi, \gamma^\vee) \, \varpi_\gamma$, and 
we can choose the coordinates $(\xi, \gamma^\vee)$ arbitrarily. Now 
given $\xi' \in \CM'$ we want the corresponding $\xi \in \CM$ to 
satisfy 
\begin{equation}
(\xi, \gamma^\vee) = \left\{\begin{array}{lr} 0 & \text{if } \gamma \in I (\lam),\\
(\xi', \gamma^\vee) & \text{if } \gamma \in J (\lam). 
\end{array}\right. \label{3.16}
\end{equation}
This makes sure that the first and the third condition in~(\ref{3.13}) 
are satisfied. 


It remains to define all $(\xi, \gamma^\vee)$ with $\gamma \notin 
\Delta_\lam$, i.e., with $(\lam, \gamma^\vee) \neq 0$, such that 
also the second condition in~(\ref{3.13}) holds. Here we require:
\begin{equation}
\hbox
{If $\gamma \in \Delta \setminus \Delta_\lam$ is linked to an element 
in~$J' (\xi') \cup \{ \alp_0\}$, then $(\xi, \gamma^\vee) > 0$.} \label{3.17}
\end{equation}

Let us show that (\ref{3.17}) implies $J (\xi) \subset J (\lam)$. Consider 
$\gamma \in J (\xi)$. The connectedness of $J (\xi)\cup\{\alpha_0\}$ implies that there exists 
a sequence $\gamma = \gamma_0, \gamma_1, \ldots, \gamma_r 
= \alp_0$ with $\gamma_i \in J (\xi)$ and $(\gamma_i, \gamma_{i+1}^\vee) 
< 0$ for all $i < r$.

If $r = 1$, then $(\xi, \gamma^\vee) = 0$ --- since $\gamma \in J (\xi)$ --- 
and (\ref{3.17}) imply $\gamma \in \Delta_\lam$, hence $\gamma \in J (\lam)$ 
by the definition of~$J (\lam)$. 

We now use induction on~$r$. So we may assume that $\gamma_i 
\in J (\lam)$ for all~$i$ with $0 < i < r$. We get then from (\ref{3.16}) that
$(\xi', \gamma_i^\vee) = (\xi, \gamma_i^\vee) = 0$ whenever
$0 < i < r$, hence $\gamma_i \in \Delta'_{\xi'}$. And since
$\{ \gamma_1, \gamma_2, \ldots, \alp_0 \}$ is connected, we get
even $\gamma_i \in J' (\xi')$, in particular $\gamma_1 \in J' (\xi')$.
Now $(\xi, \gamma^\vee) = 0$ and (\ref{3.17}) imply $\gamma \in \Delta_\lam$, 
hence $\gamma \in J (\lam)$ by the definition of~$J (\lam)$. 

\Remark: It is clear that $\CM$ is not uniquely determined by the conditions
in~(\ref{3.13}) since they only involve $\Delta_\xi$ --- which then determines 
$J (\xi)$ --- and not the precise values of the $(\xi, \gamma^\vee)$ 
with $\gamma \notin \Delta_\xi$. We shall exhibit later on a 
natural choice for~$\CM$ with nice properties. 

\medskip\noindent
{\it Examples\/}: (1) In each case $\CM'$ contains the largest short 
root~$\beta'_1$ of~$\Phi'$. Let us show that we can take the largest 
short root~$\beta_1$ of~$\Phi$ as the corresponding element in~$\CM$.
We have observed above that  $J (\beta_1) = \emptyset$. Any root 
$\gamma \in I (\lam)$ satisfies $(\alp_0, \gamma^\vee) = 0$, hence
$\gamma \in \Delta_{\beta_1}$. So there is no problem with the
first two  conditions in~(\ref{3.13}). 

Next, $\Delta_{\beta_1}$ consists exactly of the simple roots not linked 
to~$\alp_0$ in the extended Dynkin diagram, hence $J (\lam) \cap 
\Delta_{\beta_1}$ of those simple roots not linked to~$\alp_0$ in the 
extended Dynkin diagram and belonging to~$J (\lam)$. 

The same applies to  $J (\lam) \cap \Delta'_{\beta'_1}$ since 
$\Delta' \cup \{ \alp_0 \} = J (\lam) \cup \{ \alp_0, \alp' \}$ is the
extended Dynkin diagram for~$\Phi'$. So also the third condition
in~(\ref{3.13}) holds.

\smallskip\noindent
(2) \ Consider $\Phi$ of type~$F_4$ with $J (\lam) = \{ \alp_2, 
\alp_3, \alp_4 \}$. This actually implies $\lam = \frac 1 2  \varpi_1$ 
and $I (\lam) = \emptyset$. In the extended Dynkin diagram 
$\alp_0$ is linked to~$\alp_4$. This shows that $J (\lam) \cup 
\{ \alp_0 \}$ has type~$C_4$. In the extended Dynkin diagram
of  $J (\lam) \cup \{ \alp_0 \}$ the extra vertex~$\alp'$  is linked 
to~$\alp_4$. The pair $(\Delta', \Delta' \setminus \{ \alp' \})$ has 
type~$(C_4, C_3)$. If we denote by $\alp'_i$ and $\varpi'_i$ the 
simple roots and the fundamental weights for~$\Phi'$, then 
we have in standard numbering $\Delta' = \{ \alp'_1 = \alp', \, 
\alp'_2 = \alp_4, \, \alp'_3 = \alp_3, \, \alp'_4 = \alp_2 \}$ and 
by {\tt Table~2} $\CM' = \{ \varpi'_2, \, 2 \, \varpi'_1 \}$. 
Here $\varpi'_2 = \beta'_1$ in the notation from above; so we can take
$\beta_1 = \varpi_4$ as the corresponding element in~$\CM$. 
Consider next $\xi' = 2 \, \varpi'_1$. Here (\ref{3.16}) says that 
we shall take $(\xi, \alp_i^\vee) = 0$ for all $i \geq 2$. Since 
$\alp_1 \notin \Delta_\lam$ and $J' (\xi') = \{ \alp_2, 
\alp_3, \alp_4 \}$, we should require $(\xi, \alp_1^\vee) > 0$ 
according to~(\ref{3.17}). The most natural choice is therefore $\xi 
= \varpi_1$. So we end with $\CM = \{  \varpi_4, \varpi_1 \}$.


\section{The Kostant cascade}\label{KC}

We keep the assumptions on $\Phi$, $\Delta$, and~$W$ with
indecomposable~$\Phi$. 
In case all roots have the same length, we call all roots short. We 
normalise the $W$--invariant scalar product on the real span~$\Err \Phi$ 
of~$\Phi$ such that $(\alp, \alp) = 2$ for all {\bf short} roots~$\alp$. 
Then any $\lam$ in the integral span~$\Zet \Phi$ (the root lattice) 
of~$\Phi$ satisfies $(\lam, \lam) \in 2\Zet $, and $(\lam, \lam) = 2$
if and only if $\lam$ is a short root.

Our normalisation implies that $(\beta^\vee, \lam) = (\beta, \lam)$ for 
all short roots~$\beta$, but  $(\beta^\vee, \lam) = \frac 1 2  (\beta, \lam)$ 
or $(\beta^\vee, \lam) = \frac 1 3  (\beta, \lam)$ if $\beta$ is a long root
and $\Phi$ of BCF-type or of type~$G_2$.

\medskip
We now introduce the Kostant cascade~$\CK$ of~$\Phi$, actually a 
variation of the usual one. We start with the largest short root 
which we denote by~$\beta_1$. Remove all simple roots~$\alp$ 
with $(\beta_1, \alp) \neq 0$ from~$\Delta$. Decompose the remaining simple 
roots (regarded as subsets of the Dynkin diagram) into connected components 
$\Delta_1, \Delta'_1, \ldots$ \ Denote by $\beta_2, \beta'_2, \ldots$ the
largest short root in the root subsystem $\Phi \cap \Zet \Delta_1$, $\Phi 
\cap \Zet \Delta'_1$, $\ldots$ \ Remove from~$\Delta_1$ all~$\alp$ with 
$(\beta_2, \alp) \neq 0$, remove from~$\Delta'_1$ all~$\alp$ with 
$(\beta'_2, \alp) \neq 0$, and so on. Decompose the remaining simple 
roots  into connected components $\Delta_2, \Delta'_2, \ldots$ \  Set $\beta_3, 
\beta'_3, \ldots$ equal to the largest short root in $\Phi \cap \Zet \Delta_2, 
\Phi \cap \Zet \Delta'_2, \ldots\,$ respectively. Continue like this until there 
are no short simple roots left. Now the Kostant cascade~$\CK$\/ is defined 
as the set of all the largest short roots encountered in this process: 
$$\CK = \{\, \beta_1, \beta_2, \beta'_2,\ldots, \beta_3, \beta'_3, \ldots \,\}.$$ 
The construction implies that the elements in~$\CK$ are pairwise orthogonal. 

Our definition here differs in the case of two root lengths from the usual one 
where one always takes the largest root. It follows that $\{ \beta^\vee \mid 
\beta \in \CK \}$ is basically the usual Kostant cascade for the dual root 
system~$\Phi^\vee$. (There are differences in types $G_2$ and~$C_n$ with 
$n$ odd, where $\alp_n^\vee$ occurs in the usual Kostant cascade 
for~$\Phi^\vee$, while we here ignore $\alp_n$ because it is long.)

\medskip
Recall that the support of a linear combination 
of simple roots is defined as 
$$
\supp \> \left(\sum_{i=1}^n m_i \alp_i\right) = 
\left\{\, \alp_i \mid m_i \neq 0 \,\right\}.
$$ 
Recall also the usual order 
relation $\leq$ on the real span of the roots where $\lam \leq \mu$ 
if and only if $\mu - \lam = \sum_{i=1}^n m_i \alp_i$ with all $m_i 
\in \Zet$, $m_i \geq 0$. 

For any $\beta \in \CK$ set $$\Delta \, (\beta) = \supp \beta.$$
Note that $\Delta \, (\beta)$ is then also the connected subset 
of~$\Delta$ from the construction of the Kostant cascade such that 
$\beta$ is the largest short root in~$\Phi \cap \Zet \Delta \, (\beta)$. 
Note that $(\beta, \alp) \geq 0$ for all $\alp \in \Delta \, (\beta)$ since 
the largest short root is dominant. In the construction of the Kostant 
cascade one can therefore replace the condition $(\beta, \alp) \neq 0$ 
by $(\beta, \alp) > 0$.

\medskip
We call an element $\beta \in \CK$ a {\it predecessor\/} of an element 
$\beta' \in \CK$ if $\Delta \, (\beta')$ is a connected component 
of $$\Delta \, (\beta) \setminus \{\, \alp \in \Delta \, (\beta) \mid (\beta, 
\alp) > 0 \,\}.$$ One has then $\Delta \, (\beta') \subsetneq \Delta \, (\beta)$ 
and $\beta' < \beta$. (The largest short root in $\Phi \cap \Zet \, \Delta \, 
(\beta)$ is larger than the largest short root in the proper subsystem 
$\Phi \cap \Zet \, \Delta \, (\beta')$.) Each element in~$\CK$ not equal 
to~$\beta_1$ has a unique predecessor. 

We define a partial ordering $\preccurlyeq$ on~$\CK$ such that $\beta 
\preccurlyeq \beta'$ if and only if $\beta = \beta'$ or there exists a chain
$\beta = \beta^{(1)}, \beta^{(2)}, \ldots, \beta^{(s)} = \beta'$ in~$\CK$ 
such that each $\beta^{(i)}$, $1 \leq i < s$, is the predecessor 
of~$\beta^{(i+1)}$. The discussion above and the construction 
of the Kostant cascade show for all $\beta, \beta' \in \CK$ that 
\begin{equation}
\beta \preccurlyeq \beta' \iff \beta' \leq \beta \iff \Delta \, (\beta') 
\subset \Delta \, (\beta).\label{4.1}
\end{equation}
For any $\beta \in \CK$ the set of all $\beta' \in \CK$ with $\beta' \preccurlyeq 
\beta$ is totally ordered and has the form $$\beta_1 = \beta^{(1)}
\preccurlyeq  \beta^{(2)} \preccurlyeq  \cdots \preccurlyeq   \beta^{(r-1)} 
\preccurlyeq  \beta^{(r)} = \beta$$ such that each $\beta^{(i)}$, $1 \leq i < r$, 
is the predecessor of~$\beta^{(i+1)}$.


For our main goal we shall be interested in the supports of the 
difference $\beta_1 - \beta$ for $\beta  \in \CK$. One can check 
case-by-case for all $\beta \neq \beta_1$ that $\supp \, (\beta_1 - \beta)$ 
has two connected components in type~$A$ and is connected 
in the other types. We prefer to give a general argument and a version 
that works in all cases. For this we look at (our version of) 
the extended Dynkin diagram of~$\Phi$. 
Denote by $\Phi^\aff$ the affine root system with Dynkin diagram
$\Delta \cup \{ \alp_0 \}$. In the simply laced case this is the 
root system of the usual untwisted affine algebra. In the other cases 
we look at the types that Kac denotes by $D_{l+1}^{(2)}$, $A_{2l-1}^{(2)}$, 
$E_6^{(2)}$, or $D_4^{(3)}$, see \cite{Ka}, \S 4.8, Tables Aff2 and Aff3.

When working with $\Phi^\aff$ then $\alp_0$ is no longer the negative 
of the largest short root~$\beta_1$ in~$\Phi$. It is now part of a basis 
for $\Err \Phi^\aff = \Err \, \alp_0 \oplus \Err \Phi$. The basic imaginary 
root~$\delta$ in~$\Phi^\aff$ is given by $\delta = \alp_0 + \beta_1$. 

The bilinear form $(\>,\>)$ extends to a positive semidefinite bilinear 
form on~$\Err \Phi^\aff$ such that $(\alp_0, \alp_0) = 2$ and $(\delta, 
\mu) = 0$ for all $\mu \in \Err \Phi^\aff$. 

The short real roots in~$\Phi^\aff$ are characterised as the elements
$\gamma$ in the root lattice $\Zet \, (\Delta \cup \{ \alp_0 \})$ with 
$(\gamma, \gamma) = 2$. For any such root $\gamma$ also 
$\delta - \gamma$ is a short real root in~$\Phi^\aff$ as $(\delta - \gamma, 
\delta - \gamma) = (\gamma, \gamma) = 2$ and since $\delta - \gamma$
again belongs to the root lattice. This applies in particular to all 
short roots in~$\Phi$. 

\medskip

\begin{lemma} \label{L4.1} Let $\beta \in \Phi$ be a short root. Then

\Teil (a) $\supp \, (\beta_1 - \beta) \cup \{ \alp_0 \}$ is a connected subset 
of the extended Dynkin diagram, and 

\Teil (b) for each $\alp \in \Delta$ with $(\beta, \alp) > 0$ also  
$\supp \, (\beta_1 - \beta) \cup \{ \alp_0 \} \cup \{ \alp \}$ is connected. 
 
\end{lemma}

\Proof: (a) \ Since $\delta - \beta$ is a short real root in~$\Phi^\aff$,
as observed above, its support in $\Delta \cup \{ \alp_0 \}$ is connected.
 The claim follows since $\delta - \beta = \alp_0 + (\beta_1 - \beta)$
has support $\supp \, (\beta_1 - \beta) \cup \{ \alp_0 \}$. 

\smallskip\noindent
(b) \  The reflection~$s_\alp$ associated to~$\alp$ maps 
$\beta$ to a short root of the form $s_\alp \, \beta = \beta - m \, \alp$
with $m$ a positive integer. Then $\beta_1 - s_\alp \, \beta =
(\beta_1 - \beta) + m \, \alp$ satisfies $\supp \, (\beta_1 - s_\alp \, \beta) 
= \supp \, (\beta_1 -  \beta)  \cup \{ \alp \}$. Now apply~(a). 

\section{Dominant weights associated to the Kostant cascade}\label{dominant}

\bigskip
We are going to construct certain dominant weights that will play 
a crucial role in our main theorem. For any $\beta \in \CK$ set 
\begin{equation}
\xi_\beta = \sum_{\beta' \preccurlyeq \beta} \beta'. \label{5.1}
\end{equation} 
We are interested in those $\xi_\beta$ that are dominant and shall 
determine 
\begin{equation}
\Xi = \{\, \xi_\beta \mid \beta \in \CK, \, \hbox {$\xi_\beta$ 
dominant} \,\} \label{5.2}
\end{equation} 
case-by-case.  Note that $\xi_{\beta_1} = \beta_1 \in \Xi$ since the largest short root
is dominant. Also observe that Joseph \cite{J} constructs  for each $\beta \in \CK$ a dominant weight in $\beta+\sum_{\beta'\preccurlyeq \beta}\mathbf{N}\beta'$. One can check that this element coincides with $\xi_\beta$ if and only if $\xi_\beta$ is dominant.

We shall at the same time determine for each $\xi_\beta \in \Xi$ 
\begin{equation}
\supp^* \xi_\beta := \supp \, (\beta_1 - \beta) = \bigcup_{\beta' \preccurlyeq
\beta} \supp \, (\beta_1 - \beta'). \label{5.3}
\end{equation} 
In order to check the 
second equal sign, note that $\beta' \preccurlyeq \beta$ implies 
$\beta \leq \beta'$, hence $\beta_1 - \beta'  \leq \beta_1 - \beta$
and thus  $\supp \, (\beta_1 - \beta')  \subset \supp \, (\beta_1 - \beta)$. 
Note also: If $m$ is the number of~$\beta' \in \CK$ with $\beta' 
\preccurlyeq \beta$, then 
\begin{equation}
(\xi_\beta, \xi_\beta) = 2 \, m \qquad  \hbox{and} 
\qquad \supp^* \xi_\beta = \supp \, (m \, \beta_1 - \xi_\beta). \label{5.4}
\end{equation}
We have clearly $\supp^* \xi_\beta = \emptyset$ for $\beta=\beta_1$ in all cases. 

A detailed case-by-case description of $\CK$ and $\Xi$ can be found in Appendix~\ref{caseproof}.

\begin{proposition} \label{PA.1} Let $\beta \in \CK$. Then $\xi_\beta$ is 
dominant if and only if $\supp \, (\beta_1 - \beta) \neq \Delta$. 
If so, then $J (\xi_\beta) = \supp \, (\beta_1 - \beta)$.

\end{proposition}

This could be proved by inspecting all cases in Appendix~\ref{caseproof}. We prefer to give a general proof that minimises case-by-case arguments. We shall proceed via a series of lemmas. Only the first one requires 
case-by-case consideration.

\begin{lemma} \label{LA.2} Let $\beta \in \CK$ such that $\beta_1$ is 
the predecessor of~$\beta$. If $\alp \in \Delta$ with 
$(\beta, \alp) > 0$, then $\alp \notin \supp \, (\beta_1 - \beta)$.
\end{lemma}

\Proof{\it :} Since $(\beta, \alp) > 0$ implies that $\alp$ occurs in the 
support of~$\beta$, the claim is obvious if $\alp$ occurs with 
coefficient~$1$ in~$\beta_1$. This is always the case in types
$A$ and~$B$. In types $C$ and~$D$ this works for $\alp$ an 
endpoint of the Dynkin diagram. The only other $\alp$ in types 
$C$ and~$D$ is $\alp = \alp_4$ where one gets $\supp \, 
(\beta_1 - \beta) = \{ \alp_1, \alp_2, \alp_3 \}$. For the 
exceptional types look at the lists in Appendix~\ref{caseproof}.

\begin{lemma} \label{LA.3} Let $\beta, \beta' \in \CK$ such that $\beta'$ is the 
predecessor of~$\beta$. Then we have $(\beta + \beta', \alp^\vee) = 0$
for all $\alp \in \Delta \, (\beta') \setminus \Delta \, (\beta)$.
\end{lemma}

\Proof: ({\it cf.} Lemma~2.7 in \cite{J}) Recall that $\Delta \, (\beta)$ 
is a connected component of $$\Delta^0 (\beta') = \{\, \gamma \in \Delta \, 
(\beta') \mid (\beta', \gamma) = 0 \,\}.$$ So if $\alp \in \Delta \, (\beta') 
\setminus \Delta \, (\beta)$ satisfies $(\beta', \alp^\vee) = 0$, then 
$\alp$ belongs to another connected component, which then 
implies that $\alp$ is orthogonal to all roots in~$\Delta \, (\beta)$, 
hence $(\beta, \alp^\vee) = 0$ and $(\beta + \beta', \alp^\vee) = 0$.

So suppose that $(\beta', \alp^\vee) \neq 0$. Since $\beta'$ is dominant 
in~$\Phi \cap \Zet \, \Delta \, (\beta')$, we get then  $(\beta', \alp^\vee) 
> 0$, and since $\beta'$ is short, $(\beta', \alp^\vee) = 1$. 
(We could a priori have $(\beta', \alp^\vee) = 2$ and $\beta' = \alp$. But then 
$\Delta \, (\beta') = \{ \alp \}$, and $\beta'$ cannot have a successor.) 
We have to prove that $(\beta, \alp^\vee) = - (\beta', \alp^\vee) = -1$. 

Suppose first that $\Delta \, (\beta')$ is not of type~$A_n$ with $n \geq 2$. 
Then $\alp$ is uniquely determined and $\Delta \, (\beta)$ is a connected
component of $\Delta \, (\beta') \setminus \{ \alp \}$. 
Since $\beta$ belongs to $\Phi \cap \Zet \, \Delta \, (\beta')$ and since $\beta'$ 
is the only dominant short root in this subsystem, there exists $\gamma \in 
\Delta \, (\beta')$ with $(\beta, \gamma^\vee) < 0$. We have clearly $\gamma 
\notin \Delta \, (\beta)$, and $\gamma$ cannot belong to one of the other 
components of $\Delta \, (\beta') \setminus \{ \alp \}$ since the other 
components are orthogonal to~$\Delta \, (\beta)$. It follows that 
$\gamma = \alp$. And since $\beta$ is short,
we get $(\beta, \alp^\vee) = -1$ as desired. 

For $\Delta \, (\beta')$ of type $A_n$ we need $n \geq 3$ in order for a 
successor to exist. Here there are two simple roots $\alp, \alp' \in \Delta \,
(\beta')$ with $(\beta', \alp^\vee) =  1 = (\beta', {\alp'}^\vee)$ and we have
$\beta = \beta' - \alp - \alp'$ which easily yields
$(\beta, \alp^\vee) = - 1 = (\beta, {\alp'}^\vee)$.

\begin{lemma} \label{LA.4} Let $\beta \in \CK$. 

\Teil (a) If $\xi_\beta$  is dominant, then $\xi_\beta = \sum_{\alp \in \Delta \, 
(\beta)} (\beta, \alp^\vee) \, \varpi_\alp$. 

\Teil (b)  If $\xi_\beta$  is dominant,  then $\xi_{\beta'}$ is 
dominant for all $\beta' \in \CK$ with $\beta' \preccurlyeq \beta$.

\Teil (c) Let  $\beta' \in \CK$ be the predecessor of~$\beta$. 
If  $\xi_{\beta'}$ is dominant, then $\xi_\beta$  is dominant
if and only if $(\beta, \alp^\vee) = 0$ for all $\alp \in \Delta \setminus 
\Delta \, (\beta')$.

\end{lemma}

\Proof: We use induction on~$\beta_1 - \beta$. If $\beta = \beta_1$,
then (a) is clear since $\Delta \, (\beta_1) = \Delta$ while (b) and~(c) 
are empty.

So suppose that $\beta \neq \beta_1$. Let  $\beta' \in \CK$ be the 
predecessor of~$\beta$. Any $\alp \in \Delta \, (\beta)$ is 
not only orthogonal to~$\beta'$, 
but also to all elements in~$\CK$ constructed earlier. This implies 
$$(\xi_{\beta'}, \alp^\vee) = 0 \qquad  \hbox{and} 
\qquad (\xi_\beta, \alp^\vee) = (\xi_{\beta'} + \beta, \alp^\vee) = 
(\beta, \alp^\vee) \geq 0.$$ 

If $\alp \in \Delta \setminus \Delta \, (\beta)$, then $\alp 
\notin \supp \beta$, hence $(\beta, \alp^\vee) \leq 0$ and 
$$(\xi_{\beta'}, \alp^\vee) =  (\xi_\beta - \beta, \alp^\vee)
\geq (\xi_\beta, \alp^\vee).$$ This shows that 
$\xi_{\beta'}$ is dominant if $\xi_\beta$ is so. 
The more general claim in~(b) follows by induction.

Suppose now conversely that $\xi_{\beta'}$ is dominant. Since
$(\xi_\beta, \alp^\vee)  \geq 0$ for all $\alp \in \Delta \, (\beta)$, 
we see that $\xi_\beta$ is dominant if and only if $(\xi_\beta, \alp^\vee)  
\geq 0$ for all $\alp \in \Delta \setminus \Delta \, (\beta)$. 

If $\alp \in \Delta \, (\beta') \setminus \Delta \, (\beta)$, then 
we have $(\xi_{\beta'}, \alp^\vee) = (\beta', \alp^\vee)$, hence 
$$(\xi_\beta, \alp^\vee) = (\xi_{\beta'} + \beta, \alp^\vee) = 
(\beta' + \beta, \alp^\vee) = 0$$ thanks to Lemma~\ref{LA.3}.

Recall that $\lam = \sum_{\alp \in \Delta} (\lam, \alp^\vee) \, \varpi_\alp$ 
for any $\lam$. If $\alp \in \Delta \setminus \Delta \, (\beta')$, then (a) 
applied inductively to $\xi_{\beta'}$ yields $(\xi_{\beta'}, \alp^\vee) = 0$, 
hence $$(\xi_\beta, \alp^\vee) = (\xi_{\beta'} + \beta, \alp^\vee) = 
(\beta, \alp^\vee) \leq 0,$$ where the final inequality follows from 
the fact that $\alp \notin \supp \beta$. It follows that $\xi_\beta$ 
is dominant if and only if this inequality is an equality for each~$\alp$. 
Therefore (c) holds. 

Now (a) follows from the results above on all $(\xi_\beta, \alp)$.

\Remark: Note that the condition in~(c) is equivalent to
$$\{\, \alp \in \Delta \mid (\beta, \alp) < 0 \,\} = 
\{\, \alp \in \Delta \, (\beta') \mid (\beta', \alp) > 0 \,\}.
$$

\begin{lemma}\label{LA.5} Let $\beta, \beta' \in \CK$ such that $\beta'$ is the 
predecessor of~$\beta$. Suppose that $\xi' := \xi_{\beta'}$ is dominant 
and that $J (\xi') = \supp \, (\beta_1 - \beta)$. 

\Teil (a) We have either $\Delta \, (\beta) \subset J (\xi')$ 
or  $\Delta \, (\beta) \subset I (\xi')$. 
 
\Teil (b) If  $\Delta \, (\beta) \subset I (\xi')$, then $\xi_\beta$ 
is dominant. 

\Teil (c) If $\Delta \, (\beta) \subset J (\xi')$,  then $\xi_\beta$ 
is not dominant. 

\end{lemma}

\Proof: (a) We have $\Delta_{\xi'} = J (\xi') \cup I (\xi')$, hence 
$$\Delta \, (\beta) = (\Delta \, (\beta) \cap (J (\xi') \cup \{ \alp_0 \})) 
\cup (\Delta \, (\beta) \cap I (\xi')).$$ Since $\Delta \, (\beta)$ is 
connected, one of the two intersections on the right hand side 
has to be empty.

\smallskip\noindent
(b) Suppose that $\Delta \, (\beta) \subset I (\xi')$. Let $\gamma 
\in \Delta \setminus \Delta \, (\beta')$. We have $\gamma \in 
\supp \, (\beta_1 - \beta')$ by the construction of the Kostant cascade,
hence $\gamma \in J (\xi')$ by our assumption. 

We want to show $(\alp, \gamma) = 0$ for all $\alp \in \Delta \, (\beta)$, 
since then also $(\beta, \gamma) = 0$, and Lemma~\ref{LA.4}(c) implies that 
$\xi_\beta$ is dominant. 

Suppose by contradiction that there exists $\alp \in \Delta \, (\beta)$ 
with $(\alp, \gamma) \neq 0$. Then $\{ \alp, \gamma \}$ is a connected
subset of $\Delta_{\xi'} \cup \{ \alp_0 \}$ and has a non-trivial intersection
with the connected component $J (\xi') \cup \{ \alp_0 \}$, hence is
contained in $J (\xi') \cup \{ \alp_0 \}$. This yields $\alp \in J (\xi')$ 
in contradiction with $\alp \in \Delta \, (\beta) \subset I (\xi')$. 

\smallskip\noindent
(c) Suppose that $\Delta \, (\beta) \subset J (\xi')$. We want to show
that there exist $\gamma \in \Delta \setminus \Delta \, (\beta')$ and
$\alp \in \Delta \, (\beta)$ with $(\alp, \gamma) < 0$. If so, then also
$(\beta, \gamma) < 0$, hence $\xi_\beta$ is not dominant by Lemma~\ref{LA.4}(c).

Start with an arbitrary $\alp \in \Delta \, (\beta)$. Since $\Delta \, (\beta) 
\subset J (\xi')$ and  since $J (\xi') \cup \{ \alp_0 \}$ is connected, 
there exists a sequence $\alp^{(0)} = \alp_0, \alp^{(1)}, \ldots, \alp^{(r)} 
= \alp$ with $(\alp^{(i)}, \alp^{(i+1)}) < 0$ for all~$i < r$ and with 
$\alp^{(i)} \in J (\xi')$ for all $i > 0$. By changing $\alp$, we may assume 
that $\alp^{(r-1)} \notin \Delta \, (\beta)$. Note that $(\alp_0, \alp^{(1)}) < 0$ 
implies $(\beta_1, \alp^{(1)}) > 0$, hence $\alp^{(1)} \notin \Delta \, 
(\beta)$. So we have $r > 1$, hence $\alp^{(r-1)} \neq \alp_0$
and $\alp^{(r-1)} \in J (\xi') \setminus \Delta \, (\beta)$. Now 
$(\xi', \alp^{(r-1)}) = 0$ and  $(\alp^{(r-1)}, \alp^{(r)}) < 0$ imply: 
If $\alp^{(r-1)} \in \Delta \, (\beta')$, then $\alp^{(r-1)}$ belongs to the
connected component $\Delta \, (\beta)$ of $\Delta^0 (\beta')$ 
contradicting $\alp^{(r-1)} \notin \Delta \, (\beta)$. It follows that 
$\gamma \in \Delta \setminus \Delta \, (\beta')$, so we can take 
$\gamma = \alp^{(r-1)}$. 

\begin{lemma} \label{LA.6} Let $\beta \in \CK$.

\Teil (a) If $\xi_\beta$ is dominant, then $\alp \notin \supp \, (\beta_1 
- \beta)$ for all $\alp \in \Delta (\beta)$ with $(\beta, \alp) > 0$; 
we have then $J (\xi_\beta) = \supp \, (\beta_1 - \beta)$.

\Teil (b) If $\xi_\beta$ is not dominant, then $\supp \, (\beta_1 - \beta) 
= \Delta$. 

\end{lemma}

\Proof: (a) \ We use induction on $\beta_1 - \beta$. The case $\beta = 
\beta_1$ is trivial. So assume that $\beta \neq \beta_1$ and that 
$\xi_\beta$ is dominant. Denote the predecessor of~$\beta$ by~$\beta'$.
Then $\xi_{\beta'}$ is dominant by Lemma~\ref{LA.4}(b); we may assume by 
induction that $J (\xi_{\beta'}) = \supp \, (\beta_1 - \beta')$.

Now Lemma~\ref{LA.5} implies $\Delta \, (\beta) \subset I (\xi_{\beta'})$, hence
$\Delta \, (\beta) \cap \supp \, (\beta_1 - \beta') =  \Delta \, (\beta) \cap 
J (\xi_{\beta'}) = \emptyset$. It follows that $\Delta \, (\beta) \cap 
\supp \, (\beta_1 - \beta) = \Delta \, (\beta) \cap \supp \, (\beta' - \beta)$ 
and that $\supp \, (\beta_1 - \beta') \subset \Delta_{\xi_\beta}$ by Lemma~\ref{LA.4}(a).

We now apply Lemma~\ref{LA.3} 
to the root system $\Phi \cap \Zet \, \Delta \, (\beta')$. It follows that 
all $\alp \in \Delta \, (\beta)$ with $(\beta, \alp) > 0$ do not belong 
to $\supp \, (\beta' - \beta)$. It follows that $\supp \, (\beta' - \beta) 
\subset \Delta_{\xi_\beta}$, hence $\supp \, (\beta_1 - \beta) 
\subset \Delta_{\xi_\beta}$. Since $\supp \, (\beta_1 - \beta) \cup 
\{ \alp_0 \}$ is connected by Lemma~\ref{L4.1}(a), we get $\supp \, (\beta_1 
- \beta) \cup  \{ \alp_0 \} \subset J (\xi_\beta) \cup \{ \alp_0 \}$, 
i.e., $\supp \, (\beta_1 - \beta) \subset J (\xi_\beta)$. 

We want to prove equality. Suppose that we have a root $\alp \in 
J (\xi_\beta) \setminus \supp \, (\beta_1 - \beta)$. Since  $\Delta \setminus 
\Delta \, (\beta) \subset \supp \, (\beta_1 - \beta) $,
we get $\alp \in \Delta \, (\beta)$. We have then $(\beta_1, \alp) = 0$ 
as $\alp \in \Delta \, (\beta)$ and $(\beta, \alp) = 0$ as $\alp \in 
J (\xi_\beta)$, hence $(\beta_1 - \beta, \alp) = 0$. As $\alp \notin 
\supp \, (\beta_1 - \beta)$ this implies $(\alp', \alp) = 0$ for all 
$\alp'  \in \supp \, (\beta_1 - \beta)$. 

As $J (\xi_\beta) \cup \{ \alp_0 \}$ is connected and contains~$\alp$,
there is a sequence $\gamma_1 = \alp, \gamma_2, \ldots, 
\gamma_{r-1},$ $\gamma_r = \alp_0$ with $\gamma_i \in J (\xi_\beta)$ 
for all $i < r$ and $(\gamma_i, \gamma_{i+1}) < 0$ for all $i < r$. 
As $(\beta_1, \gamma_{r-1}) = - (\alp_0, \gamma_{r-1}) > 0$, 
we cannot have $\gamma_{r-1} \in \Delta \, (\beta)$, hence 
$\gamma_{r-1} \in \supp \, (\beta_1 - \beta)$. 

Let $i > 1$ be minimal for $\gamma_i \in \supp \, (\beta_1 - \beta)$. 
Then $\gamma_{i-1} \in J (\xi_\beta) \setminus \supp \, (\beta_1 - \beta)$ 
implies as above $(\gamma_{i-1}, \gamma_i) = 0$ 
contradicting our choice of the sequence. So there is no $\alp$ as above.

\smallskip\noindent
(b) Let $\beta \in \CK$ be the smallest element with $\beta' \preccurlyeq 
\beta$ such that $\xi_{\beta'}$ is not dominant. Since $\beta_1 = 
\xi_{\beta_1}$ is dominant, we have $\beta\ \neq \beta_1$, so 
$\beta'$ has a predecessor~$\beta''$. Then $\xi_{\beta''}$ is 
dominant by the minimality of~$\beta'$. Thanks to~(a) we can now
apply Lemma~\ref{LA.5} and get $\Delta \, (\beta') \subset J (\xi_{\beta''}) = 
\supp \, (\beta_1 - \beta'')$, hence $$\Delta \, (\beta) \subset 
\Delta \, (\beta') \subset \supp \, (\beta_1 - \beta'') \subset \supp \, 
(\beta_1 - \beta).$$
Since any $\alp \in \Delta \setminus \Delta \, (\beta)$ clearly 
belongs to~$\supp \, (\beta_1 - \beta)$. 

\Remark: Comparing both parts of the lemma we see
$$\hbox {$\xi_\beta$ dominant} \iff \supp \, (\beta_1 - \beta) 
\neq \Delta.$$ Together with the first part of the lemma this 
shows that we have proved Proposition~\ref{PA.1}.

\section{Proof of the main theorem}\label{proofMT}

We start with another look at the minuscule case.

\begin{lemma} \label{L5.2} Let $\alp \in \Delta$ such that $\varpi_\alp$ 
is minuscule. We have then for any $\beta \in \CK$ that 
\begin{equation}
\alp \in \Delta \, (\beta) \iff \alp \notin \supp \, (\beta_1 - \beta).
\label{5.5}
\end{equation} 
The set $\CM$ from {\tt Table 2} for~$\alp$ is equal 
to the set of all $\xi_\beta$ with $\alp \in \Delta \, (\beta)$. 
\end{lemma}

\Proof: Recall that a fundamental weight~$\varpi_\alp$
is minuscule if and only if $1 = (\varpi_\alp, \beta_1^\vee) = 
(\varpi_\alp, \beta_1)$. The condition $(\varpi_\alp, \beta_1^\vee) 
= 1$ implies that $\alp$ occurs with coefficient~$1$ when we write 
$\beta_1$ as a linear combination of the simple roots. This fact
yields the equivalence in~(\ref{5.5}). The claim on the set from
{\tt Table~2} follows by inspection. 

\Remark: The lemma implies that all $\xi_\beta$ with $\alp \in 
\Delta (\beta)$ are dominant. 

\medskip
Fix now $\lam \in \Ahh$ with $(\lam, \beta_1^\vee) = 1$. 
Let us write $$\lam = \sum_{i=1}^n r_i \, \varpi_i.$$ So we have 
$r_i \geq 0$ for all~$i$ and $\Delta_\lam = \{\, \alpha_i \mid 1 \leq i \leq n, \> 
r_i = 0 \,\}$. Since $\beta \leq \beta_1$ for all $\beta \in \CK$, 
we get $0 \leq (\beta, \lam) \leq (\beta_1, \lam) = 1$  and $$(\beta, \lam) 
= 1 \iff \hbox {$r_i = 0$ for all $\alp_i \in \supp (\beta_1 - \beta)$} 
\iff \supp \, (\beta_1 - \beta) \subset \Delta_\lam.$$
Set now 
$$
\Xi (\lam) := \{\, \xi \in \Xi \mid \supp^* \xi \subset \Delta_\lam \,\}.
\label{5.6}
$$ 
If $\lam = \varpi_\alp$ is a minuscule fundamental 
weight, then $\Delta_{\varpi_\alp} = \Delta \setminus \{ \alp\}$, 
so $\supp^* \xi \subset \Delta_\lam$ is equivalent to $\alp \notin 
\supp^* \xi$. Therefore Lemma~\ref{L5.2} shows that $\Xi (\varpi_\alp)$ 
is the set~$\CM$ from {\tt Table~2}.

\begin{lemma} \label{L5.3} We have $$\Xi (\lam) = \{\, \xi \in \Xi \mid 
(\xi, \xi) = 2 \, (\xi, \lam) \,\}.$$
\end{lemma}

\Proof: Suppose $\xi = \xi_\beta$ with $\beta \in \CK$. Set $m$ 
equal to the number of~$\beta' \in \CK$ with $\beta' \preccurlyeq \beta$. 
We have then $(\xi, \xi) = 2 \, m$ and $$(\xi, \lam) = \sum_{\beta' 
 \preccurlyeq \beta} \, (\beta', \lam) \leq m$$ with equality if and 
 only if $\supp \, (\beta_1 - \beta') \subset \Delta_\lam$ for all $\beta' 
 \preccurlyeq \beta$, hence by~(\ref{5.3})  if and  only if $\supp^* \xi 
 \subset \Delta_\lam$.

\begin{lemma} \label{L5.4} If $\xi \in \Xi (\lam)$, then $J (\xi) \subset J (\lam)$ 
and $I (\lam) \subset \Delta_\xi$.
\end{lemma}

\Proof: We have $J (\xi) \cup \{ \alp_0 \} = \supp^* \xi \cup \{ \alp_0 \} 
\subset \Delta_\lam \cup \{ \alp_0 \}$. Since the left hand side is 
connected, it is contained in the connected component $J (\lam) 
\cup \{ \alp_0 \}$ of the right hand side containing~$\alp_0$. 
This implies $J (\xi) \subset J (\lam)$. 

Consider now $\alp \in \Delta_\lam$ with $\alp \notin \Delta_\xi$, 
hence $(\xi, \alp) > 0$. We have to show that $\alp \in J (\lam)$. Let 
$\beta \in \CK$ with $\xi = \xi_\beta$, hence with $J (\xi) = \supp \, (\beta_1 
- \beta)$ by Lemma~\ref{LA.6}. We have then $(\beta, \alp) = (\xi_\beta, \alp) > 0$ 
by Lemma \ref{LA.4}(a). We get in particular $\alp \in 
\Delta \, (\beta)$, so Lemma~\ref{L4.1}(b) shows that $$\supp \, (\beta_1 - \beta) 
\cup \{ \alp_0 \} \cup \{ \alp \}  = J (\xi) \cup \{ \alp_0 \} \cup \{ \alp \}$$ is 
connected. As it is a subset of $\Delta_\lam \cup \{ \alp_0 \}$, it is 
already contained in $J (\lam) \cup \{ \alp_0 \}$. This yields
$\alp \in J (\lam)$. 

We are now ready to prove our main theorem, slightly reformulated in a more compact form as below:
\medskip\noindent
\begin{thm}[reformulated] {\it We have} \label{Tmainr}
\begin{equation}
P (W_\lam^{\rm aff} /W_\lam) = 1 
+ \sum_{\xi \in \Xi (\lam)} \, t^{\ell_\beta} \, P (W_\lam / 
W_{\lam, \xi}). \label{5.7}
\end{equation}
\end{thm}

\Proof: If $\lam$ is a minuscule fundamental weight, then the claim 
follows from Lemma~\ref{L5.2}. In general, we have to show that we can take
$\CM = \Xi (\lam)$ in Proposition~\ref{P3.6}. Recall the construction there.
We consider the extended Dynkin diagram associated to $J (\lam) \cup 
\{ \alp_0 \}$, constructed by adding a vertex~$\alp'$. We consider the 
root system~$\Phi'$ with Dynkin diagram $\Delta'  = J (\lam) \cup 
\{ \alp' \}$. Set $\CM'$ equal to the set of weights for~$\Phi'$ associated
to $\Phi'$ and~$\alp'$ in {\tt Table~2}. By Lemma~\ref{L3.7} it now 
suffices to find a bijection $\Xi (\lam) \isto \CM'$, $\xi \mapsto \xi'$,
satisfying~(\ref{3.13}). Note that the first two conditions in~(\ref{3.13}) hold 
thanks to Lemma~\ref{L5.4}. So we only have to check that $J (\lam) \cap 
\Delta_\xi = J (\lam) \cap \Delta'_{\xi'}$ or, equivalently, that 
$J (\lam) \setminus \Delta_\xi = J (\lam) \setminus \Delta'_{\xi'}$.  

We proceed case-by-case looking at the distinct possibilities for 
the type of the pair $(J (\lam) \cup \{ \alp_0 \}, J (\lam))$, equal to 
the type of one of the pairs $(\Phi', \Delta' \setminus \{ \alp' \})$. 
We shall denote the simple roots in~$\Delta'$ and the corresponding
fundamental weights by $\alp'_i$ and~$\varpi'_i$ respectively, 
using standard numbering.

\medskip\noindent
{\bf Type $\bf (A_{k+1}, A_k)$}: This possibility corresponds to the 
subcase $j = 1$ of the type-A case in {\tt Table 1}. (It includes the 
case $J (\lam) = \emptyset$ for $k = 0$.) In this case the largest 
short root~$\beta'$ of~$\Phi'$ is the only element of~$\CM'$. 
The largest short root~$\beta_1 = \xi_1$ belongs to (any!) 
$\Xi (\lam)$ and the pair $(\beta_1, \beta'_1)$ satisfies~(\ref{3.13}) 
according to the first example following Proposition~\ref{P3.6}.

So we only have to show that $\beta_1$ is the only element 
in~$\Xi (\lam)$. If $\xi \in \Xi (\lam)$, then the connected subset 
$J (\xi) \cup \{ \alp_0 \}$ is contained in $J (\lam) \cup \{ \alp_0 \}$.
Since $\alp_0$ is an end-vertex of $J (\lam) \cup \{ \alp_0 \}$, 
the pair $(J (\xi) \cup \{ \alp_0 \}, J (\xi))$ has to have type 
$(A_{s+1}, A_s)$ for some $s \leq k$. However {\tt Table~3} 
over all elements in~$\Xi$ shows that $\xi = \beta_1$ is the 
only element in~$\Xi$ where this type occurs.
  
\medskip\noindent
{\bf Type $\bf (A_{j+k+1},A_j \times A_k)$} with $k \geq j > 0$: 
Since $J (\lam)$ is not connected, this possibility can occur 
only for $\Phi$ of type~$A_n$ for some~$n$. We may assume 
up to symmetry that $$J (\lam) = \{ \alp_1, \alp_2, \ldots, \alp_j, 
\alp_{n+1-k}, \ldots, \alp_{n-1}, \alp_n\}$$ and get $$\Delta' = 
\{ \alp'_1 = \alp_1, \alp'_2 = \alp_2, \ldots, \alp'_j = \alp_j, 
\alp'_{j+1} = \alp', \alp'_{j+2} = \alp_{n+1-k}, \ldots, \alp'_{j+k+1} = 
\alp_n\}.$$ We get now from {\tt Table~2} and from a look 
at all $\supp^* \xi$ that $$\CM' = \{\, \varpi'_i + \varpi'_{j+k+2-i} 
\mid 1 \leq i \leq j+1 \,\} \quad \hbox {and} \quad \Xi (\lam) = 
\{\, \varpi_i + \varpi_{n+1-i} \mid 1 \leq i \leq j +1\,\}.$$
Then we have a bijection given by $\xi_i := \varpi_i + \varpi_{n+1-i} 
\mapsto  \xi'_i :=  \varpi'_i + \varpi'_{j+k+2-i}$ and it has the right
property since $$J (\lam) \setminus \Delta_{\xi_i} = \{\, \alp_i, \, 
\alp_{n+1-i} \,\} =  \{\, \alp'_i, \ \alp'_{j+k+2-i} \,\} = J (\lam) \setminus 
\Delta'_{\xi'_i}$$ if $i < j+1$. For $i = j+1 < k+1$ both sets 
are equal to $\{ \alp_{n+1-j} \}$, and for $i = j+1 = k+1$ both sets 
are empty. 
 
\medskip\noindent
{\bf Type $\bf (B_{k+1}, A_k)$} with $k > 0$: This possibility can occur 
only for $\Phi$ of type~$B_n$ with $n \geq k+1$, and we get $J (\lam) 
= \{ \alp_1, \alp_2, \ldots, \alp_k \}$. We may assume that $k < n-1$ 
since otherwise $\lam$ is the minuscule fundamental weight~$\varpi_n$.  
A look at all $\supp^* \xi$ shows that $\Xi (\lam) = \{\, \varpi_i \mid  
1 \leq i \leq k+1\,\}$. We get $\Delta' = \{ \alp'_1 = \alp_1, \alp'_2 
= \alp_2, \ldots, \alp'_k = \alp_k, \alp'_{k+1} = \alp' \}$ and now 
{\tt Table 2} yields $\CM' = \{\, \varpi'_i \mid 1 \leq i \leq k \} \cup 
\{ 2 \, \varpi'_{k+1} \}$. Now the obvious bijection with $\varpi_i 
\mapsto \varpi'_i$ for $i \leq k$ and $\varpi_{k+1} \mapsto 
2 \, \varpi'_{k+1}$ works.

\medskip\noindent
{\bf Type $\bf (D_{k+1}, A_k)$} with $k \geq 3$: This possibility can occur 
only for $\Phi$ of type $C_n$ or~$D_n$ with $n \geq k+1$, and we get 
$J (\lam) = \{ \alp_1, \alp_2, \ldots, \alp_k \}$ (up to symmetry in 
type~$D_n$). Let us assume $k \leq n-2$ in case $\Phi$ has type~$D_n$
so to exclude the case that $\lam$ is the minuscule fundamental 
weight~$\varpi_n$. A look at all $\supp^* \xi$ shows that $\Xi (\lam) 
= \{\, \varpi_{2i} \mid  1 \leq i \leq (k+1)/2\,\}$ where we have to replace
$\varpi_{n-1}$ in type~$D_n$ by $\varpi_{n-1} + \varpi_n$. 

We have $\Delta' = \{ \alp'_1 = \alp_1, \alp'_2 = \alp_2, \ldots, \alp'_k 
= \alp_k, \alp'_{k+1} = \alp' \}$. (Note that $\alp'$ is linked to~$\alp_{k-1}$ 
in the diagram.) We get from {\tt Table~2} that $\CM' = \{\, \varpi'_{2i} 
\mid 1 \leq i \leq (k+1)/2 \}$ where we have to replace $\varpi'_k$ 
by $\varpi'_k + \varpi'_{k+1}$ for $k$ even, and $\varpi'_{k+1}$ 
by $2 \, \varpi'_{k+1}$ for $k$ odd. Now the map $\varpi_{2i} 
\mapsto \varpi'_{2i}$ with the obvious modification works. 

\medskip\noindent
{\bf Types  $\bf (C_{k+1}, C_k)$ and $\bf (D_{k+1}, D_k)$}:  
Here $\CM'$ has order~$2$ and contains besides $\beta'_1 
= \varpi'_2$ the element~$2 \, \varpi'_1$. We have here $\alp'_1 
= \alp'$, hence $\Delta'_{2 \varpi'_1} = J (\lam)$. So we have to 
show that also $\Xi (\lam)$ consists of two elements: $\beta_1 
= \xi_1$ corresponding to $\beta'_1$, and another one, say 
$\xi$, that has to satisfy $J (\lam) \subset \Delta_\xi$. 

Except for $\Phi$ of type $C_{k+1}$ or~$D_{k+1}$ with $\lam$ 
a minuscule fundamental weight, there are four possibilities 
for~$\Phi$: 

\smallskip
\begin{description}
\item{$\bullet$} $\Phi$ has type $F_4$ with 
$J (\lam) = \{ \alp_2, \alp_3, \alp_4 \}$,

\item{$\bullet$} $\Phi$ has type~$E_6$ with $J (\lam) = \{ \alp_2, 
\alp_3, \alp_4, \alp_5 \}$,

\item{$\bullet$} $\Phi$ has type~$E_7$ with $J (\lam) = \{ \alp_1, \alp_2, 
\alp_3, \alp_4, \alp_5 \}$,

\item{$\bullet$} $\Phi$ has type~$E_8$ with $J (\lam) = \{ \alp_2, 
\alp_3, \alp_4, \alp_5, \alp_6, \alp_7, \alp_8 \}$. 
\end{description}

\smallskip\noindent
In each case a look at all $\supp^* \xi$ shows that $\Xi (\lam) = 
\{ \xi, \xi_2 \}$ with $\xi_2 = \varpi_1 / \varpi_1 + \varpi_6 / \varpi_6 
/ \varpi_1$, hence with $J (\lam) \subset \Delta_{\xi_2}$. 

\medskip\noindent
{\bf Types  $\bf (E_6, D_5)$ and $\bf (E_7, E_6)$}:  One checks 
that these possibilities occur only for $\Phi$ of type $E_6$ or 
$E_7$ respectively with $\lam$ a minuscule fundamental weight.
So we do not have to work.

\bigskip

\appendix
\section{Tables for $\CK$ and $\Xi$}\label{caseproof}
\bigskip
Below we usually write $\xi_i = \xi_{\beta_i}$, $\xi'_i = \xi_{\beta'_i}$, and so on.

\medskip\noindent
{\bf Type A$_{\bf n}$}: \ Here we have 
$$\CK = \{\, \beta_i = \eps_i - 
\eps_{n+2-i} = \alp_i + \alp_{i+1} + \cdots + \alp_{n+1-i} \mid 1 \leq i 
< (n+2)/2\,\}.$$ 
The set is totally ordered: We have $\beta_1
\preccurlyeq  \beta_2 \preccurlyeq  \cdots \preccurlyeq  \beta_m$ 
where $m = [n/2]$, and get 
$$\Xi = \{\, \xi_i  = \beta_1 + \beta_2 + \cdots 
+ \beta_i = \varpi_i + \varpi_{n+1-i}\}.$$ 
(For $n$ odd
this means $\xi_{(n+1)/2} = 2 \, \varpi_{(n+1)/2}$.) The support
is given by $$\supp^* \xi_i = \{\, \alp_1, \alp_2, \ldots, \alp_{i-1},
\alp_{n+2-i}, \ldots, \alp_{n-1}, \alp_n \,\}.$$

\medskip\noindent
{\bf Type B$_{\bf n}$}, $n \geq 2$: Here we have 
$$
\CK = \{\, \beta_i = \eps_i 
= \alp_i + \alp_{i+1} + \cdots + \alp_n \mid 1 \leq i \leq n \,\}.
$$ 
The set is totally ordered: We have $\beta_1
\preccurlyeq  \beta_2 \preccurlyeq  \cdots \preccurlyeq  \beta_n$
and get 
$$
\Xi = \{\, \xi_i = \beta_1 + \beta_2 + \cdots + \beta_i 
\mid 1 \leq i \leq n \,\},  \qquad \xi_i = \left\{\begin{array}{lr} \varpi_{i} & \text{ if } 1 \leq i < n,\\
2 \, \varpi_n & \text{ if } i = n.
\end{array} \right.
$$ 
The support is given by 
$$
\supp^* \xi_i = 
 \{ \alp_1, \alp_2, \ldots, \alp_{i-1} \}.
 $$ 
 

\medskip\noindent
{\bf Type C$_{\bf n}$}, $n \geq 3$: Here $\CK$ is the union of 
$$\{\,  \beta_i 
= \eps_{2i-1} + \eps_{2i} = \alp_{2i-1} + 2 \, \alp_{2i} + \cdots + 
2 \, \alp_{n-1}  + \alp_n \mid 1 \leq i \leq n/2\,\}$$ 
and 
$$
\{\, \beta'_i 
= \eps_{2i-3} - \eps_{2i-2} = \alp_{2i-3} \mid 2 \leq i \leq (n+2)/2 \,\}.
$$
 
We have $\beta_1 \preccurlyeq  \beta_2 \preccurlyeq  \cdots \preccurlyeq  
\beta_m$ with $m = [n/2]$ and $\beta_i \preccurlyeq  \beta'_{i+1}$
for all~$i$, $1 \leq i \leq n/2$. We get 
$$
\Xi = \{\, \xi_i = \beta_1 + 
\beta_2 + \cdots + \beta_i = \varpi_{2i} \mid 1 \leq i \leq n/2 \,\} 
\cup \{ \xi'_2 = \beta_1 + \beta'_2 = 2 \, \varpi_1 \}.
$$ 
The support
is given by (for $i > 1$) 
$$\supp^* \xi_i =   \{\, \alp_1, \alp_2, \ldots, 
\alp_{2i-1} \,\} \qquad  \hbox{and} \qquad \supp^* \xi'_2 =  \{\, \alp_2, 
\alp_3, \ldots, \alp_n \,\}.
$$


\medskip\noindent
{\bf Type D$_{\bf n}$}, $n \geq 4$: Here $\CK$ is the union of  
$$
\{\, \beta_i 
= \eps_{2i-1} + \eps_{2i} \mid 1 \leq i \leq  n/2\,\}, \; \beta_i 
= \left\{ \begin{array}{ll} \alp_{2i-1} + 2 \sum_{j=2i}^{n-2} \alp_j + 
\alp_{n-1} + \alp_n & \text{ if } 2i \leq n-1\\
 \alp_n & \text{ if } 2i=n,
 \end{array}\right.
 $$ 
and 
$$
\{\, \beta'_i = \eps_{2i-3} - \eps_{2i-2} = \alp_{2i-3} \mid 2 \leq  i 
\leq (n+1)/2 \,\}
$$ 
and, if $n = 2 \, m$ is even, 
$$
\beta''_m = \eps_{n-1} - \eps_n 
= \alp_{n-1}.
$$ 
We have $\beta_1 \preccurlyeq  \beta_2 \preccurlyeq  \cdots 
\preccurlyeq  \beta_m$ with $m = [n/2]$ and $\beta_i \preccurlyeq  
\beta'_{i+1}$ for all~$i$, $1 \leq i \leq [(n-1)/2]$ and for even $n = 2 \, m$
also $\beta_{m-1} \preccurlyeq  \beta''_m$. Now $\Xi$ is the union of 
$$
\{\, \xi_i = \beta_1 + \beta_2 + \cdots + \beta_i \mid 1 \leq i \leq n/2 \,\},
\qquad \xi_i = \left\{ \begin{array}{ll}  \varpi_{2i} & 
\text{if } 2i < n-1, \\
 \varpi_{n-1} + \varpi_n& \text{if } 2i = n-1,\\
  2 \varpi_n& \text{if }   2i = n,
  \end{array}\right. 
$$ 
together with 
$$
\xi'_2 = \beta_1 
+ \beta'_2 =  2 \, \varpi_1
$$ 
and for even $n = 2m$ also $$\xi''_m 
= \beta_1 + \beta_2 + \cdots + \beta_{m-1} + \beta''_m = 2 \, \varpi_{n-1}.$$ 
The support is given by  (for $i > 1$) 
$$
\supp^* \xi_i =   \{\, \alp_1, \alp_2, 
\ldots, \alp_{2i-1} \,\} \qquad  \hbox{and} \qquad \supp^* \xi'_2 =  \{\, \alp_2, 
\alp_3, \ldots, \alp_n \,\}
$$  
and (in case $n = 2 \, m$ is even) $\supp^* 
\xi''_m = \{\, \alp_1, \alp_2, \ldots, \alp_{n-2}, \alp_n \,\}$.


\bigbreak
In the exceptional cases we are going to list only those
$\beta \in \CK$ where $\supp \, (\beta_1 - \beta)$ is not the 
full set of all simple roots. The other ones are not of interest 
for our applications; one can also check that the corresponding
$\xi_\beta$ are not dominant, cf.\ also Lemma \ref{LA.6}(b) . 
The remaining~$\beta \in \CK$ are linearly ordered under~$\preccurlyeq$. 

\smallskip\noindent
{\bf Type E$_{\bf 6}$}: \  Besides $\beta_1 =  \alp_1 + 2 \, \alp_2 + 2 \, \alp_3 
+ 3 \, \alp_4 + 2 \, \alp_5 + \alp_6 = \varpi_2 =\xi_1$ we have 
$\beta_2 = \alp_1 + \alp_3 + \alp_4 + \alp_5 + \alp_6 = 
\varpi_1 + \varpi_6 - \varpi_2$ with corresponding element in~$\Xi$ given
by  $\xi_2 = \beta_1 + \beta_2 = \varpi_1 + \varpi_6$. 
We get $$\supp^* \xi_2 = \{\, \alp_2, \alp_3 , \alp_4, \alp_5 \,\}.$$ 

\smallskip\noindent
{\bf Type E$_{\bf 7}$}: \ Besides $\beta_1 =  2 \, \alp_1 + 2 \, \alp_2 + 
3 \, \alp_3 + 4 \, \alp_4 + 3 \, \alp_5 + 2 \, \alp_6  + \alp_7 = \varpi_1 =\xi_1$ 
we have $\beta_2 = \alp_2 + \alp_3 + 2 \,  \alp_4 + 
2 \, \alp_5 + 2 \, \alp_6 + \alp_7 = \varpi_6 - \varpi_1$ and 
$\beta_3 = \alp_7 = 2 \, \varpi_7 - \varpi_6$. The corresponding 
elements in~$\Xi$ are $\xi_2 = \beta_1 + \beta_2 = \varpi_6$ 
and $\xi_3 = \beta_1 + \beta_2 + \beta_3 = 2 \, \varpi_7$.  
We get $$\supp^* \xi_2 = \{\, \alp_1, \alp_2, \alp_3 , \alp_4, \alp_5 \,\}, \qquad
\supp^* \xi_3 = \{\, \alp_1, \alp_2, \alp_3 , \alp_4, \alp_5, \alp_6 \,\}.$$ 

\smallskip\noindent
{\bf Type E$_{\bf 8}$}: \ Besides $\beta_1 =  2 \, \alp_1 + 3 \, \alp_2 + 
4 \, \alp_3 + 6 \, \alp_4 + 5 \, \alp_5 + 4 \, \alp_6  + 3 \, \alp_7 + 2\, \alp_8 
= \varpi_8 =\xi_1$ we have 
$\beta_2 = 2 \, \alp_1 + 2 \, \alp_2 + 3 \,\alp_3 + 
4 \, \alp_4 + 3 \, \alp_5 + 2 \, \alp_6 + \alp_7 = \varpi_1 - \varpi_8$. The
corresponding element in~$\Xi$ is $\xi_2 = \beta_1 + \beta_2 
= \varpi_1$. We get 
$$
\supp^* \xi_2 =  \{\, \alp_2, \alp_3 , \alp_4, \alp_5, \alp_6,
\alp_7, \alp_8 \,\}.
$$

\smallskip\noindent
{\bf Type F$_{\bf 4}$}: \ Besides $\beta_1 = \alp_1 + 2 \, \alp_2 + 3 \, \alp_3 
+ 2 \, \alp_4 = \varpi_4 = \xi_1$ we have $\beta_2 = \alp_1 + \alp_2 + \alp_3 = 
\varpi_1 - \varpi_4$. The corresponding element in~$\Xi$ is 
$\xi_2 = \beta_1 + \beta_2 =  \varpi_1$. We get $$\supp^* \xi_2 =  \{\, \alp_2, 
\alp_3 , \alp_4 \,\}.$$

\smallskip\noindent
{\bf Type G$_{\bf 2}$}: \ Here only $\beta_1 = 2 \alp_1 + \alp_2 = 
\varpi_1 = \xi_1$ occurs.

\bigskip

\bigskip 
If we compare these lists with {\tt Table~3} above, then we will 
see that the weights in that table are just the elements in~$\Xi$ 
for a given type. And a comparison with the results on
$\supp^* \xi$ yields:

\begin{lemma}
We have that $J(\xi) = \text{supp}^*\xi$  for each $\xi\in \Xi$.
\end{lemma}

\section{Alternative proof for the main theorem}
\bigskip

First, we are going to look in more detail at the set of all~$\beta \in \CK$
with $\supp \, (\beta_1 - \beta) \subset I$ for some proper subset~$I$
of~$\Delta$. The results will be applied to $I = J (\lam)$.

\begin{lemma} \label{LB.1} Let $I \subset \Delta$ be a proper subset. 
Let $\beta \in \CK$ with $\supp \, (\beta_1 - \beta) \subset I$. 
Then each $\beta' \in \CK$ with $\beta' \preccurlyeq \beta$
also satisfies $\supp \, (\beta_1 - \beta') \subset I$. There 
exists at most one root $\beta'' \in \CK$ with $\beta$ as 
predecessor and $\supp \, (\beta_1 - \beta'') \subset I$.

\end{lemma}

\Proof: For the first claim, recall that $\supp \, (\beta_1 - \beta') 
\subset \supp \, (\beta_1 - \beta)$ for $\beta'  \preccurlyeq \beta$,
see right after~(\ref{5.3}). 

If $\beta'' \in \CK$ has $\beta$ as predecessor, then there exists 
a connected component~$J$ of $\Delta (\beta) \cap \Delta_\beta$
with $\Delta (\beta'') = \supp \beta'' = J$. If $\supp \, (\beta_1 - \beta'') 
\subset I$, then we get $\Delta \setminus J \subset \supp \, (\beta_1 
- \beta'') \subset I$.

If also $\gamma \in \CK$ has $\beta$ as predecessor, then 
$\Delta (\gamma) \cap J = \emptyset$ by the construction of~$\CK$.
It follows that $J \subset \supp \, (\beta_1 - \gamma)$. If now also
$\supp \, (\beta_1 - \gamma) \subset I$, then we would get 
$\Delta = J \cup (\Delta \setminus J) \subset I$ --- a contradiction 
since we assume $I$ to be proper. So~$\beta''$ is unique (if 
it exists). 

\Remark: This shows that the roots $\beta \in \CK$ with $\supp \, (\beta_1 
- \beta) \subset I$ form a complete interval $\beta_1 \preccurlyeq \beta_2 
\preccurlyeq \cdots \preccurlyeq \beta_r$ in~$\CK$ with respect 
to~$\preccurlyeq$. 

\begin{lemma} \label{LB.2} Let $I \subset \Delta$ be a proper subset, let $\alp
\in \Delta \setminus I$. 

\Teil (a) If $\beta \in \Phi$ is a root with $\beta \leq \beta_1$ and
$\supp \, (\beta_1 - \beta) \subset I$, then $(\beta, \alp) \geq 
(\beta_1, \alp) \geq 0$.

\Teil (b) Suppose that $\beta \in \CK$ with $\supp \, (\beta_1 - \beta) 
\subset I$. If $\alp$ is not linked to any element of~$I \cup \{ \alp_0 \}$,
then $(\xi_\beta, \alp) = 0$. If $\alp$ is linked to some  element 
of~$\supp \, (\beta_1 - \beta) \cup \{ \alp_0 \}$, then $(\xi_\beta, \alp) > 0$.

\end{lemma}

\Proof: (a) Since $\alp \notin \supp \, (\beta_1 - \beta)$ and since $\beta_1
- \beta$ is a linear combination of the simple roots with 
non-negative coefficients, we have $(\beta_1 - \beta, \alp) \leq 0$.
This implies the claim since $\beta_1$ is dominant.

\smallskip\noindent
(b) Note that $\alp \notin I$, hence $\alp \notin \supp \, (\beta_1 - \beta)$, 
implies $\alp \in \supp \, \beta$. Since $\xi_\beta$ is dominant by 
Lemma~\ref{LA.6}, we get thus $(\xi_\beta, \alp) = (\beta, \alp)$. 

If $\alp$ is not linked to any $\gamma \in I \cup \{ \alp_0 \}$, then we 
have on the one hand $(\alp_0, \alp) = 0$, hence $(\beta_1, \alp) = 0$. 
On the other hand we get $(\gamma, \alp) = 0$ for all $\gamma \in 
\supp \, (\beta_1 - \beta) \subset I$. It follows that  $(\beta_1 - \beta, 
\alp) = 0$, hence also $0 = (\beta, \alp) = (\xi_\beta, \alp)$.

\medskip
Suppose now that $\alp$ is linked to some to some $\gamma \in 
\supp \, (\beta_1 - \beta) \cup \{ \alp_0 \}$. 
In case $\gamma = \alp_0$ we get $(\alp_0, \gamma) < 0$, 
hence $(\beta_1, \alp) > 0$. Now (a) yields $0 < (\beta, \alp) = 
(\xi_\beta, \alp)$. 

If $\gamma \in \supp \, (\beta_1 - \beta)$, then $(\gamma, \alp) < 0$ 
implies  $(\beta_1 - \beta, \alp) < 0$, hence $(\beta, \alp) 
> 0$ by~(a).

\bigskip\noindent
{\bf Set-up}. 
Consider the following situation: We have two indecomposable 
finite root systems $\Phi_1$ and~$\Phi_2$ with bases $\Delta_1$ 
and~$\Delta_2$. Denote by $\beta_1$ resp.\ $\beta'_1$ the largest 
short root and use it to construct the affine root system $\Phi_1^\aff$ 
resp.~$\Phi_2^\aff$. Write $\alp_0$ resp.~$\alp'_0$ for the extra 
basis element for $\Phi_1^\aff$ resp.~$\Phi_2^\aff$, and $\delta =
\alp_0 + \beta_1$ and $\delta' = \alp'_0 + \beta'_1$ for the basic 
imaginary roots.

Let now $I_i \subset \Delta_i$ be {\it proper\/} subsets for $i = 1, 2$ 
such that $I_1 \cup \{ \alp_0 \}$ and $I_2 \cup \{ \alp'_0 \}$ are 
connected. Suppose that we have a bijection $$f \colon I_1 \cup \{ \alp_0 \}
\isto I_2 \cup \{ \alp'_0 \} \qquad \hbox {with $f (\alp_0) = \alp'_0$}$$
that induces an isomorphism of Dynkin diagrams. Extend $f$ to a linear 
isomorphism $\Err (I_1 \cup \{ \alp_0 \}) \isto \Err (I_2 \cup \{ \alp'_0 \})$. 
The assumption involving the Dynkin diagram implies for the corresponding 
simple reflections $$s_{f (\alp)} \, f (\beta) = f (s_\alp \, \beta) \qquad 
\hbox {for all $\alp, \beta \in I_1 \cup \{ \alp_0 \}$.}$$ Since $\Phi_I = 
W_I I$ for any proper subset $I$ of the basis of one of the $\Phi_i^\aff$,
this shows that $f$ maps $\Phi_{I_1 \cup \{ \alp_0 \}}$ bijectively 
onto $\Phi_{I_2 \cup \{ \alp'_0 \}}$. That assumption implies also 
that $(f (\lam), f (\mu)) = (\lam, \mu)$ first for all $\lam, \mu \in 
I_1 \cup \{ \alp_0 \}$ and then for all  $\lam, \mu \in \Err \, 
(I_1 \cup \{ \alp_0 \})$. (Start with the fact that $(\alp_0, \alp_0) 
= 2 = (\alp'_0, \alp'_0)$ and use the connectedness of~$I_1 \cup 
\{ \alp_0 \}$.) 

\medskip 
In order to simplify notation, let us now identify $I_1 \cup \{ \alp_0 \}$ 
and $I_2 \cup \{ \alp'_0 \}$ via~$f$. So we write $I = I_1 = I_2$ 
and $\alp'_0 = \alp_0$. We regard $\Err (I \cup \{ \alp_0 \})$ as 
a subspace both of $\Err \, \Phi_1^\aff$ and of  $\Err \, \Phi_2^\aff$. 
The root subsystem $\Phi_{I \cup \{ \alp_0 \}}$ generated by
$I \cup \{ \alp_0\}$ is the same in~$\Phi_1^\aff$ and in~$\Phi_2^\aff$.

Set $\Psi_1$ equal to the set of all short roots~$\beta \in \Phi_1$
with $\supp \, (\beta_1 - \beta) \subset I$, set $\Psi_2$ equal to the set 
of all short roots~$\beta \in \Phi_2$ with $\supp \, (\beta'_1 - \beta) 
\subset I$. Note that $\Psi_1$ and~$\Psi_2$ consist of positive 
roots since (e.g.) $\supp (\beta_1 - \beta) = \Delta_1$ for any negative 
root~$\beta \in \Phi_1$.

\begin{lemma} \label{LB.3} There exists a bijection $\tau \colon \Psi_1 \isto \Psi_2$
such that $$\beta_1- \beta = \beta'_1 - \tau (\beta) \eqno(1)$$ for all 
$\beta \in \Psi_1$. We have $(\tau (\beta), \alp) = (\beta, \alp)$ 
for all $\alp \in I \cup \{ \alp_0 \}$ and  $(\tau (\beta), \tau (\beta')) 
= (\beta, \beta')$ for all $\beta, \beta' \in \Psi_1$.  

\end{lemma}

\Proof: Let $\beta \in \Psi_1$. Write $\beta_1 - \beta = \eta \in \Zet I$.
We observed before Lemma~\ref{L4.1} that $\gamma: = \delta - \beta$ is a short 
real root in~$\Phi_1^\aff$. We have $\gamma = \alp_0 + \beta'_1 - \beta
 = \alp_0 + \eta$, hence $\gamma \in \Phi_{I \cup \{ \alp_0 \}}$. So 
$\gamma$ is also a short real root in~$\Phi_2^\aff$. Therefore 
$\delta' - \gamma$ is another short real root in~$\Phi_2^\aff$. 
Define then~$\tau$ by $\tau (\beta) = \delta' - \gamma = \beta'  
- \eta$. We see that $\tau (\beta) \in \Phi_2$ and that (1) holds. 
And (1) implies $\supp \, (\beta' - \tau (\beta)) = \supp \, (\beta_1 
- \beta) \subset I$, hence that $\tau (\beta) \in \Psi_2$.

We get for all $\beta \in \Psi_1$ and $\alp \in I \cup \{ \alp_0 \}$ 
$$(\beta, \alp) = - (\delta - \beta, \alp) = - (\delta' - \tau (\beta), \alp) 
= ( \tau (\beta), \alp)$$ and for all $\beta, \beta' \in \Psi_1$ 
$$(\beta, \beta') = (\delta - \beta, \delta - \beta') = (\delta' - \tau (\beta), 
\delta' - \tau (\beta')) = (\tau (\beta), \tau (\beta')).$$ 

Since the set-up is symmetric in $\Delta_1$ and~$\Delta_2$ 
we get similarly a map $\tau' \colon \Psi_2 \to \Psi_1$ with 
$\beta'_1 - \gamma = \beta - \tau' (\gamma)$ for all $\gamma 
\in \Psi_2$. Then $\tau$ and~$\tau'$ are clearly inverse to each 
other, hence bijections.

\bigskip
Let us denote the Kostant cascade in~$\Phi_i$  by~$\CK_i$, $i = 1, 2$.

\begin{proposition} \label{PB.4} The bijection $\tau$ maps the set of all 
$\beta \in \CK_1$ with $\supp \, (\beta_1 - \beta) \subset I$ 
onto the set of all $\gamma \in \CK_2$ iwith $\supp \, (\beta'_1 
- \gamma) \subset I$.

\end{proposition}

\Proof: Let $\beta \in \CK_1$ with $\supp \, (\beta_1 - \beta) \subset I$.
We want to show that $\tau (\beta) \in \CK_2$. 
We use induction on the number of elements $\beta'$ in the 
Kostant cascade with $\beta' \preccurlyeq \beta$. The induction 
starts since $\tau (\beta_1) = \beta'_1$ belongs to~$\CK_2$.

Suppose now that $\beta \neq \beta_1$. Then $\beta$ has a 
predecessor~$\beta'$ that by Lemma~\ref{LB.1} satisfies 
the same assumptions. Write $\gamma = \tau (\beta)$ and $\gamma' 
= \tau (\beta')$. We may assume by induction that $\gamma'$ belongs 
to~$\CK_2$. 

Set $\eta = \beta_1 - \beta$ and $\eta' = \beta_1 - \beta'$. Now 
$\beta' \geq \beta$ implies $\eta' \leq \eta$, hence $\gamma' = 
\beta'_1 - \eta' \geq \beta'_1 - \eta = \gamma$. This implies 
in particular that $S := \supp \gamma \subset S' := \supp \gamma'$. 

Since $\gamma'$ belongs to the Kostant cascade, it is the largest
short root with support~$S'$ and thus satisfies $(\gamma', \alp) \geq 0$ 
for all $\alp \in S'$, hence also for all~$\alp \in S$. On the other hand, 
we have $(\beta', \beta) = 0$ since these roots are two distinct 
members of a Kostant cascade. Now Lemma~\ref{LB.3} implies $(\gamma', 
\gamma) = 0$.  Write $\gamma = \sum_{\alp \in S} m_\alp \alp$ 
with integers $m_\alp > 0$. Then $0 = (\gamma', \gamma) = 
\sum_{\alp \in S} m_\alp (\gamma', \alp)$ implies by the observation 
above that $(\gamma', \alp) = 0$ for all $\alp \in S$. We get thus 
$$S \subset S' \cap (\Delta_2)_{\gamma'}.\eqno(2)$$

Since $\supp \, (\beta'_1 - \gamma) = \supp \, \eta \subset I$, Lemma~\ref{LB.2}(a) 
implies $$(\gamma, \alp) \geq (\beta'_1, \alp) \geq 0 \qquad \hbox 
{for all $\alp \in  \Delta_2 \setminus I$.} \eqno(3)$$ 
We claim next that $$\hbox {$\alp \in \Delta_2$ with $(\gamma, \alp) 
< 0$} \; \Longrightarrow (\gamma', \alp) > 0. \eqno(4)$$ To start 
with $(\gamma, \alp) < 0$ implies $\alp \in I$ by~(3), 
hence $(\beta, \alp) < 0$ by Lemma~\ref{LB.3}. Since $\xi_\beta$ is dominant, 
Lemma~\ref{LA.4} (applied to~$\Delta_1$) shows that $\alp \in \supp (\beta')$. 
We get now from Lemma~\ref{LA.3} and Lemma~\ref{LB.3} $$(\gamma', \alp) = 
(\beta', \alp) = - (\beta, \alp) > 0$$ as claimed.

Let us show next that $S$ is a connected component of $S' \cap 
(\Delta_2)_{\gamma'}$.  If not, then there exist $\alp \in (S' \cap 
(\Delta_2)_{\gamma'}) \setminus S$ and $\alp' \in S$ 
with $(\alp, \alp') < 0$. This implies $(\gamma, \alp) < 0$, hence 
$(\gamma', \alp) > 0$ by~(4), contradicting $\alp \in (\Delta_2)_{\gamma'}$. 

We claim now that $\gamma$ is the largest short root in~$\Phi_{S}$,
or, equivalently, that $(\gamma, \alp) \geq 0$ for all $\alp \in S$. 
But this  follows  from~(4) since $S \subset (\Delta_2)_{\gamma'}$. 

Now $\gamma' \in \CK_2$ implies $\gamma \in \CK_2$ by the 
construction of the Kostant cascade; furthermore $\gamma'$ is
the predecessor of~$\gamma$. 

\smallskip
We have thus shown that $\tau (\CK_1 \cap \Psi_1) \subset \CK_2 
\cap \Psi_2$. Since our set-up is symmetric in $\Phi_1$ and~$\Phi_2$ 
and since $\tau^{-1}$ is defined symmetrically, we get now equality: 
$\tau (\CK_1 \cap \Psi_1) = \CK_2 \cap \Psi_2$.

\Remark: Recall that the roots $\beta \in \CK_1$ with $\supp \, (\beta_1 
- \beta) \subset I$ form a complete interval $\beta_1 \preccurlyeq \beta_2 
\preccurlyeq \cdots \preccurlyeq \beta_r$ in~$\CK_1$ with respect 
to~$\preccurlyeq$. The proof shows more precisely that $\tau (\beta_1) 
\preccurlyeq \tau (\beta_2) \preccurlyeq \cdots \preccurlyeq \tau (\beta_r)$ 
is the complete interval of all roots $\gamma \in \CK_2$ with 
$\supp \, (\beta'_1 - \gamma) \subset I$.

\bigskip

Now we return to the set-up of Proposition~\ref{P3.6}. There we had fixed $\lam \in 
\Ahh$ with $(\lam, \beta_1^\vee) = 1$. We consider there 
a root system~$\Phi'$ with basis $\Delta' = J (\lam) \cup \{ \alp' \}$.
The construction there shows that the assumptions of the 
set-up above before Lemma~\ref{LB.3} are satisfied with $\Phi_1 = 
\Phi$ and $\Phi_2 = \Phi'$ and $I = J (\lam)$. 

Since $\alp'$  corresponds to a  minuscule fundamental weight
for the root system~$\Phi'$, it occurs with coefficient~$1$ 
in~$\beta'_1$. So we have by Lemma~\ref{L5.2} for any $\beta'$ 
in the Kostant cascade for~$\Phi'$ that $$\alp' \in \supp \beta' 
\iff \supp \, (\beta'_1 - \beta') \subset J (\lam).$$ 

Denote by $\beta'_1 \preccurlyeq \beta'_2 \preccurlyeq \cdots 
\preccurlyeq \beta'_r$ the elements in the Kostant cascade for~$\Phi'$ 
containing $\alp'$ in its support. In particular, $\beta'_1$ is 
the largest short root in~$\Phi'$. Set $\xi'_j = \sum_{i=1}^j \beta'_i$ 
for each~$j$. So the elements $\xi'_j$ are dominant and form 
the set~$\CM'$ in the proof of Proposition~\ref{P3.6}.

Set $\beta_i = \tau^{-1} (\beta'_i)$ for all~$i$. Proposition~\ref{PB.4}
tells us that $\beta_1 \preccurlyeq \beta_2 \preccurlyeq \cdots 
\preccurlyeq \beta_r$ the elements~$\beta$ in the Kostant cascade 
for~$\Phi$ with $\supp \, (\beta_1 - \beta) \subset J (\lam)$. 
Let us write $\xi_i = \xi_{\beta_i}$ for $1 \leq i \leq r$. We have then
$$\Xi (\lam) = \{\, \xi_i \mid 1 \leq i \leq r \,\}.$$

\begin{lemma} \label{LB.5} We can take $\CM = \Xi (\lam)$ in Lemma~\ref{L3.7} 
with the bijection $\CM \to \CM'$ given by $\xi_i \mapsto \xi'_i$.

\end{lemma}

\Proof: We have $J (\xi_i) \subset J (\lam)$ and $I (\lam) \subset
\Delta_{\xi_i}$ for all~$i$ thanks to Lemma~\ref{L5.4}. So it remains to 
check for all~$i$ that $$J (\lam) \cap \Delta_{\xi_i} = 
J (\lam) \cap \Delta_{\xi'_i},$$ i.e., that for all $\alp \in J (\lam)$
$$(\xi_i, \alp) = 0 \iff (\xi'_i, \alp) = 0.$$ But we have even stronger
$$(\xi_i, \alp) = (\xi'_i, \alp) \qquad \hbox {for all $\alp \in J (\lam)$}$$
as Lemma~\ref{LB.3} yields $(\beta_j, \alp) = (\beta'_j, \alp)$ for all~$j$. 

\medskip\noindent

Now we are ready to give another proof of the main theorem:

{\it Proof of Theorem~\ref{Tmain}}: If $\lam$ is a minuscule fundamental weight, then the claim 
follows from Lemma~\ref{L5.2}. In general, combine Lemma~\ref{LB.5} and 
Proposition~\ref{P3.6}.

\end{document}